\documentclass[11pt]{amsart}
\usepackage{amsmath, amssymb, amsfonts, amsthm}
\usepackage{mathrsfs}
\usepackage{hyperref}
\usepackage{url}
\usepackage{enumitem}
\usepackage{geometry}
\usepackage{cite}
\title{Spectral--Geometric Deformations of Function Algebras on Manifolds}

\theoremstyle{plain}
\newtheorem{theorem}{Theorem}[section]
\newtheorem{lemma}[theorem]{Lemma}
\newtheorem{proposition}[theorem]{Proposition}
\newtheorem{corollary}[theorem]{Corollary}

\theoremstyle{definition}
\newtheorem{definition}[theorem]{Definition}
\newtheorem{assumption}[theorem]{Assumption}
\newtheorem{example}[theorem]{Example}

\theoremstyle{remark}
\newtheorem{remark}[theorem]{Remark}

\newcommand{\Spec}{\mathrm{Spec}}
\newcommand{\id}{\mathrm{id}}
\newcommand{\Hom}{\mathrm{Hom}}
\newcommand{\bbC}{\mathbb{C}}
\newcommand{\bbR}{\mathbb{R}}
\newcommand{\bbZ}{\mathbb{Z}}
\newcommand{\bbT}{\mathbb{T}}
\newcommand{\cC}{\mathcal{C}}

\newcommand{\Aut}{Aut}
\DeclareMathOperator{\ev}{ev}
\DeclareMathOperator{\suppspec}{supp_{\Delta}}

\author{Amandip Sangha}
\address{The Climate and Environmental Research Institute NILU}
\email{asan@nilu.no}

\begin{document}

\begin{abstract}
We introduce an intrinsic deformation of the algebra of smooth functions on a compact Riemannian manifold using only the Laplace spectral decomposition.  The construction twists the canonical multiplication--projection channels by unimodular phases, producing a well-defined bilinear product on the finite spectral core with values in $L^2(M)$.  We give a simple condition for compatibility with complex conjugation and isolate a Sobolev boundedness hypothesis under which the product extends to a Sobolev algebra and admits iteration; in that setting, associativity is equivalent to an explicit identity for the twisted spectral channels.  We analyze gauge and coboundary aspects for scalar twists and obtain rigidity statements in the action-free regime.  We also compare with classical strict deformation frameworks arising from actions of locally compact abelian groups---Rieffel's deformation for $\mathbb{R}^d$-actions, Connes--Landi's torus isospectral deformations, and Kasprzak's cocycle deformation via Landstad theory---showing that, when the relevant abelian group action has a discrete spectral decomposition (in particular, in the compact abelian/periodic case where the algebra decomposes into homogeneous subspaces indexed by characters of the acting group), their deformed products are recovered uniformly as refined instances of our channel twist.  Finally, we formulate a grading-based obstruction and classification for graded scalar twists.
\end{abstract}

\maketitle

\setcounter{tocdepth}{1}
\tableofcontents

\section{Introduction}

Deformation procedures for algebras of functions are often driven by extra structure: a group action, a foliation, a Poisson bracket, or a quantization prescription tied to symmetries.  This paper takes a different starting point.  On any compact Riemannian manifold, the Laplace--Beltrami operator provides a canonical orthogonal spectral decomposition, and that decomposition canonically organizes pointwise multiplication into ``channels'' obtained by multiplying and then projecting to a fixed eigenspace.

Let $P_\lambda$ be the orthogonal projection onto the Laplace eigenspace $E_\lambda$.  Multiplying two functions and projecting defines the canonical multiplication channels
\[
m^\nu_{\lambda\mu}(f,g):=P_\nu\bigl((P_\lambda f)(P_\mu g)\bigr).
\]
A basic technical point is that, even when $f$ and $g$ are finite spectral sums, the product typically has infinitely many nonzero spectral components.  Accordingly, any ``spectral deformation'' must be formulated so that infinite output support is harmless.  We fix unimodular scalars $\omega^\nu_{\lambda\mu}\in\bbT$ and twist each channel by $\omega^\nu_{\lambda\mu}$.  Summing the twisted outputs in the square-integrable sense yields a bilinear deformed product
\[
\star_\omega:\ E_{\mathrm{fin}}\times E_{\mathrm{fin}}\longrightarrow L^2(M),
\]
defined intrinsically from the Laplace projections and requiring no group action or auxiliary input.

For arbitrary twisting data we prove: (i) $\star_\omega$ is always well-defined as an $L^2$-convergent spectral series on $E_{\mathrm{fin}}$; (ii) a simple symmetry of $\omega$ ensures compatibility with complex conjugation; and (iii) there is a natural gauge action by diagonal spectral unitaries, yielding a large intrinsic family of strictly associative (indeed commutative) twists obtained by conjugating pointwise multiplication.  These statements are genuinely symmetry-free and depend only on the Laplace decomposition.

To discuss associativity beyond the core one needs a function space closed under $\star_\omega$.  Rather than hiding this analytic issue, we separate it cleanly by isolating a Sobolev boundedness hypothesis: under this hypothesis $\star_\omega$ extends continuously to a Sobolev algebra.  In that setting, associativity becomes a concrete, checkable condition: it is equivalent to an explicit channel identity on eigenvectors.  This reduces the algebraic constraint to an intrinsic statement about how spectral multiplication channels compose.

Although our construction does not assume symmetries, it is compatible with symmetry-driven strict deformation frameworks \emph{when they are present}.  Classical strict deformations are built from actions of locally compact abelian groups (Rieffel: $\bbR^d$ \cite{Rieffel}; Connes--Landi: $\bbT^d$ \cite{ConnesLandi}; Kasprzak: general locally compact abelian $G$ via Landstad theory and a cocycle on $\widehat G$ \cite{Kasprzak}).  When the relevant abelian action has a discrete spectral decomposition---for instance, for compact abelian actions (or for periodic instances of $\bbR^d$-actions that factor through a torus) where the algebra decomposes into homogeneous subspaces indexed by characters of the acting group---and isometric so that it commutes with $\Delta$---one obtains a discrete homogeneous decomposition (a grading) that refines the Laplace eigenspaces. In that setting, the classical deformed products admit a uniform homogeneous-core description by cocycle factors, and we show that these are recovered by an appropriate refinement of our channel twist.  We also make precise the structural slogan that ``grading is the essence'' of cocycle twisting: an action supplies a grading on a suitable homogeneous core, cocycle twisting depends only on the graded pieces, and different actions yield the same twisted product precisely when they induce the same grading (up to relabeling).  A C*-level corollary records that topological gradings are implemented by the dual compact abelian action, so graded cocycle twisting in the C*-setting is automatically action-based in this dual sense.

Within the graded cocycle subclass we give a cohomological obstruction/classification picture, including a rank-one obstruction for noncommutativity (Corollary~\ref{cor:rank-one}).  Together with the gauge-triviality results for intrinsic scalar twists, this clarifies both the scope of the scalar theory and the source of rigidity.

The present paper is intentionally minimal and scalar-valued.  While it provides a robust intrinsic definition, a sharp associativity criterion (under a stated analytic hypothesis), and a clear gauge and obstruction theory, producing genuinely new intrinsic noncommutative strictly associative examples appears to require richer channel data.  This motivates the development of non-scalar (matrix-valued) deformation mechanisms based on multiplicity and tensor-category methods, which we mention in the Outlook.

\section{Spectral preliminaries}

Let $(M,g)$ be a compact Riemannian manifold without boundary, and let $\Delta$ be the Laplace--Beltrami operator on complex-valued functions.
Its spectrum $\Spec(\Delta)\subset[0,\infty)$ is discrete \cite{Chavel1984}, each eigenspace
\[ E_\lambda:=\ker(\Delta-\lambda)\subset C^\infty(M)\cap L^2(M) \]
finite-dimensional, and
\begin{equation}\label{eq:L2-decomp}
L^2(M)=\bigoplus_{\lambda\in\Spec(\Delta)} E_\lambda
\end{equation}
as an orthogonal direct sum. Denote by $P_\lambda:L^2(M)\to E_\lambda$ the orthogonal projection.

\begin{definition}[Finite spectral core]
Define the algebraic spectral core
\[ E_{\mathrm{fin}}:=\bigoplus_{\lambda\in\Spec(\Delta)}^{\mathrm{alg}} E_\lambda, \]
consisting of finite spectral sums. For $f\in E_{\mathrm{fin}}$ its \emph{spectral support}
$\suppspec(f):=\{\lambda\in\Spec(\Delta):P_\lambda f\neq 0\}$ is finite.
\end{definition}

\begin{remark}\label{rem:Efin-not-algebra}
In general $E_{\mathrm{fin}}$ is not closed under pointwise multiplication (nor under our yet to be defined deformed product $\star_\omega$),
so we view $\star_\omega$ (Definition~\ref{def:star}) on $E_{\mathrm{fin}}$ as a bilinear map $E_{\mathrm{fin}}\times E_{\mathrm{fin}}\to L^2(M)$.
Algebra-level statements (associativity, isomorphisms, $*$-algebra structure) are made on Sobolev algebras $H^s(M)$ with
$s>\frac{\dim M}{2}$ once the relevant extensions exist. When the same formula extends boundedly on $H^s(M)$ for all $s>\frac{\dim M}{2}$
(as in the gauge family of \S\ref{sec:gauge}), it restricts to the Fr\'echet algebra $C^\infty(M)=\bigcap_s H^s(M)$.
\end{remark}

\begin{lemma}[Parseval in $L^2$, \cite{Chavel1984}]\label{lem:parseval-L2}
For every $f\in L^2(M)$ one has $f=\sum_{\lambda}P_\lambda f$ with convergence in $L^2$, and
\[ \|f\|_{L^2}^2=\sum_{\lambda\in\Spec(\Delta)}\|P_\lambda f\|_{L^2}^2. \]
\end{lemma}

\begin{proof}
This is the standard orthogonal decomposition for the self-adjoint operator $\Delta$ with compact resolvent.
\end{proof}

\subsection*{Sobolev spaces}
For $s\in\bbR$ define
\[ \|f\|_{H^s}^2:=\sum_{\lambda\in\Spec(\Delta)} (1+\lambda)^s\,\|P_\lambda f\|_{L^2}^2, \]
and let $H^s(M)$ be the completion of $E_{\mathrm{fin}}$ under $\|\cdot\|_{H^s}$ \cite{Hebey1996}.

\begin{lemma}[Parseval in $H^s$ \cite{Taylor2011}]\label{lem:parseval-Hs}
For every $s\in\bbR$ and $f\in H^s(M)$ one has $f=\sum_{\lambda}P_\lambda f$ with convergence in $H^s$, and
\[ \|f\|_{H^s}^2=\sum_{\lambda\in\Spec(\Delta)} (1+\lambda)^s\,\|P_\lambda f\|_{L^2}^2. \]
In particular $E_{\mathrm{fin}}$ is dense in $H^s(M)$.
\end{lemma}

\begin{proof}
By construction, $H^s(M)$ is the Hilbert $\ell^2$-direct sum of the weighted spaces $(1+\lambda)^{s/2}E_\lambda$.
\end{proof}

\begin{remark}[Pointwise multiplication on $H^s$]
If $s>\frac{\dim M}{2}$ then $H^s(M)$ is a Banach algebra under pointwise multiplication.
We will use this repeatedly, notably in \S\ref{sec:gauge} and \S\ref{sec:obstruction}.
\end{remark}

\section{Multiplication channels and evaluation}

For $\lambda,\mu,\nu\in\Spec(\Delta)$ define the canonical multiplication component
\[ m^{\nu}_{\lambda\mu}(f,g):=P_\nu(fg) \in E_\nu,\qquad f\in E_\lambda,\ g\in E_\mu. \]
Since $E_\lambda,E_\mu\subset C^\infty(M)$, the product $fg$ is smooth (hence in $L^2$), so $P_\nu(fg)$ is well-defined.
Equivalently, $m^{\nu}_{\lambda\mu}$ determines a linear map
\[ m^{\nu}_{\lambda\mu}\in\Hom(E_\lambda\otimes E_\mu,E_\nu). \]

\begin{definition}[Minimal channel line]
Define the (possibly zero) \emph{channel line}
\[ \cC^{\nu}_{\lambda\mu}:=\mathrm{span}_{\bbC}\{m^{\nu}_{\lambda\mu}\}\ \subset\ \Hom(E_\lambda\otimes E_\mu,E_\nu). \]
\end{definition}

\begin{definition}[Evaluation map]\label{def:ev}
Define the evaluation map
\[ \ev^{\nu}_{\lambda\mu}:\cC^{\nu}_{\lambda\mu}\otimes(E_\lambda\otimes E_\mu)\to E_\nu,
\qquad m\otimes w\mapsto m(w). \]
Here $m\in\cC^{\nu}_{\lambda\mu}$ is a \emph{linear map} $E_\lambda\otimes E_\mu\to E_\nu$.
\end{definition}

\begin{lemma}\label{lem:ev-well}
The map $\ev^{\nu}_{\lambda\mu}$ is well-defined and linear.
\end{lemma}

\begin{proof}
Since $\cC^{\nu}_{\lambda\mu}\subset\Hom(E_\lambda\otimes E_\mu,E_\nu)$, every $m\in\cC^{\nu}_{\lambda\mu}$ is a linear map
and $m(w)$ is defined. Linearity in each factor implies linearity on the tensor product.
\end{proof}

\section{Scalar twisting coefficients and the deformed product}

\subsection{Twisting data}

\begin{definition}[Scalar twisting system]
A \emph{scalar twisting system} is a family
\[ \omega= (\omega^{\nu}_{\lambda\mu})_{\lambda,\mu,\nu\in\Spec(\Delta)},\qquad \omega^{\nu}_{\lambda\mu}\in\bbT. \]
\end{definition}

Given a scalar twisting system $\omega$, define the twisted channel maps
\[ \widetilde m^{\nu}_{\lambda\mu}:=\omega^{\nu}_{\lambda\mu}\,m^{\nu}_{\lambda\mu}\in\Hom(E_\lambda\otimes E_\mu,E_\nu). \]

\subsection{Definition on $E_{\mathrm{fin}}$ and $L^2$-well-definedness}

\begin{definition}[Deformed product on the spectral core]\label{def:star}
For $f,g\in E_{\mathrm{fin}}$, define
\[ f\star_\omega g := \sum_{\nu\in\Spec(\Delta)}\sum_{\lambda,\mu\in\Spec(\Delta)} \widetilde m^{\nu}_{\lambda\mu}(P_\lambda f,P_\mu g), \]
where the sum over $(\lambda,\mu)$ is finite. We let 
\[ 
P_\nu(f\star_\omega g) := \sum_{\lambda,\mu\in\Spec(\Delta)} \widetilde m^{\nu}_{\lambda\mu}(P_\lambda f,P_\mu g),
\]
and then we can write
\[
f\star_\omega g = \sum_{\nu\in\Spec(\Delta)} P_\nu(f\star_\omega g).
\]
\end{definition}

\begin{proposition}[$L^2$-well-definedness]\label{prop:L2-well}
For all $f,g\in E_{\mathrm{fin}}$ the series $\sum_{\nu}P_\nu(f\star_\omega g)$ converges in $L^2(M)$.
Thus $\star_\omega:E_{\mathrm{fin}}\times E_{\mathrm{fin}}\to L^2(M)$ is a well-defined bilinear map.
\end{proposition}

\begin{proof}
Write $f=\sum_{\lambda\in F} f_\lambda$ and $g=\sum_{\mu\in G} g_\mu$ with finite sets $F,G\subset\Spec(\Delta)$.
For fixed $(\lambda,\mu)\in F\times G$, set $h_{\lambda\mu}:=f_\lambda g_\mu\in L^2(M)$.
By Lemma~\ref{lem:parseval-L2},
\[ \sum_{\nu}\|P_\nu h_{\lambda\mu}\|_{L^2}^2=\|h_{\lambda\mu}\|_{L^2}^2. \]
Multiplying each $P_\nu h_{\lambda\mu}$ by the unimodular scalar $\omega^{\nu}_{\lambda\mu}$ preserves the $\ell^2$-norm.
Hence $\sum_{\nu} \omega^{\nu}_{\lambda\mu}P_\nu h_{\lambda\mu}$ converges in $L^2(M)$.
Finally,
\[ P_\nu(f\star_\omega g)=\sum_{(\lambda,\mu)\in F\times G} \omega^{\nu}_{\lambda\mu}P_\nu h_{\lambda\mu} \]
is a finite sum in $E_\nu$, so $(P_\nu(f\star_\omega g))_{\nu}$ is a finite sum of $\ell^2$-sequences.
Thus $\sum_\nu P_\nu(f\star_\omega g)$ converges in $L^2(M)$.
Bilinearity is immediate.
\end{proof}

\begin{remark}[Untwisted case]
If $\omega^{\nu}_{\lambda\mu}\equiv 1$, then $P_\nu(f\star_\omega g)=P_\nu(fg)$ for all $\nu$, hence $f\star_\omega g=fg$ in $L^2$.
\end{remark}

\section{Involution structure on the core}

Let $c(f)=\overline f$ be complex conjugation.
Since $\Delta$ has real coefficients, conjugation preserves each eigenspace and commutes with $P_\lambda$.

\begin{lemma}[Conjugation commutes with $P_\lambda$]\label{lem:conj-proj}
For every $\lambda\in\Spec(\Delta)$ and $f\in L^2(M)$, $P_\lambda(\overline f)=\overline{P_\lambda f}$.
\end{lemma}

\begin{proof}
Conjugation commutes with $\Delta$ and preserves orthogonality in $L^2$.
Since $P_\lambda$ is the orthogonal projection onto $E_\lambda$, it commutes with conjugation.
\end{proof}

\begin{definition}[$*$-compatible twisting system]\label{def:star-compat}
A twisting system $\omega$ is \emph{$*$-compatible} if
\[ \omega^{\nu}_{\mu\lambda}=\overline{\omega^{\nu}_{\lambda\mu}}\qquad\text{for all }\lambda,\mu,\nu. \]
\end{definition}

\begin{proposition}[$*$-identity on $E_{\mathrm{fin}}$]\label{prop:star-id}
If $\omega$ is $*$-compatible, then for all $f,g\in E_{\mathrm{fin}}$ one has in $L^2(M)$
\[ \overline{f\star_\omega g}=\overline g\star_\omega\overline f. \]
\end{proposition}

\begin{proof}
By bilinearity it suffices to treat $f\in E_\lambda$ and $g\in E_\mu$.
Using Lemma~\ref{lem:conj-proj} and commutativity of pointwise multiplication,
\[ P_\nu(\overline{f\star_\omega g})=\overline{P_\nu(f\star_\omega g)}=
\overline{\omega^{\nu}_{\lambda\mu}}\,P_\nu(\overline f\,\overline g)=\omega^{\nu}_{\mu\lambda}P_\nu(\overline g\,\overline f)=P_\nu(\overline g\star_\omega\overline f). \]
Equality of all spectral components implies equality in $L^2$ by Lemma~\ref{lem:parseval-L2}.
\end{proof}

\section{An intrinsic associative family: spectral gauge twists}\label{sec:gauge}

Although an arbitrary twisting system $\omega$ need not produce an associative product,
there is a natural intrinsic class of strictly associative twists obtained by conjugating the pointwise product by
unitaries diagonal in the Laplace decomposition. These twists are generally not $*$-compatible unless one imposes an
additional reality condition; see Theorem~\ref{thm:gauge-trivial}.

\begin{definition}[Spectral gauge unitaries]
Let $\varphi:\Spec(\Delta)\to\bbR$.
Define a unitary $U_\varphi:L^2(M)\to L^2(M)$ by
\[ U_\varphi\big|_{E_\lambda}=e^{i\varphi(\lambda)}\id_{E_\lambda}. \]
\end{definition}

\begin{definition}[Gauge twisting system]
Given $\varphi$, define
\begin{equation}\label{eq:gauge-u}
 \omega^{\nu}_{\lambda\mu}(\varphi):=\exp\big(i(\varphi(\lambda)+\varphi(\mu)-\varphi(\nu))\big)\in\bbT.
\end{equation}
\end{definition}

\begin{theorem}[Spectral gauge twists are scalar channel twists (gauge-triviality)]\label{thm:gauge-trivial}
Define a product on $E_{\mathrm{fin}}$ by conjugating pointwise multiplication,
\[
f\star_\varphi g \ :=\ U_\varphi^{-1}\big((U_\varphi f)\,(U_\varphi g)\big),
\qquad f,g\in E_{\mathrm{fin}}.
\]
Then for all $f,g\in E_{\mathrm{fin}}$ one has the identity in $L^2(M)$
\begin{equation}\label{eq:gauge-conjugation}
f\star_{\omega(\varphi)} g \;=\; f\star_\varphi g \;=\; U_\varphi^{-1}\big((U_\varphi f)\,(U_\varphi g)\big).
\end{equation}
In particular, applying $U_\varphi$ to \eqref{eq:gauge-conjugation} yields the core conjugacy identity
\[
U_\varphi\big(f\star_{\omega(\varphi)} g\big)=(U_\varphi f)(U_\varphi g)
\qquad (f,g\in E_{\mathrm{fin}}),
\]
i.e.\ $U_\varphi$ conjugates the twisted bilinear map on $E_{\mathrm{fin}}$ (with values in $L^2(M)$)
to pointwise multiplication.

Moreover, for every $s>\frac{\dim M}{2}$ the product extends to a bounded bilinear map
\[
\star_{\omega(\varphi)}:H^s(M)\times H^s(M)\to H^s(M),
\]
given by the same conjugation formula \eqref{eq:gauge-conjugation}.
For each such $s$, $U_\varphi$ implements an (isometric) algebra isomorphism
\[
U_\varphi:\ (H^s(M),\star_{\omega(\varphi)})\xrightarrow{\ \cong\ }(H^s(M),\cdot),
\]
so $\star_{\omega(\varphi)}$ is strictly associative (indeed commutative) on $H^s(M)$ and hence restricts
to a Fr\'echet algebra product on $C^\infty(M)=\bigcap_s H^s(M)$.

Finally, $\star_{\omega(\varphi)}$ is $*$-compatible (with respect to complex conjugation) whenever $\omega(\varphi)$ is
$*$-compatible in the sense of Definition~\ref{def:star-compat}; for instance this holds if $\varphi(\lambda)\in\pi\bbZ$
for all $\lambda$.
\end{theorem}

\begin{proof}
Let $f=\sum_{\lambda\in F} f_\lambda$ and $g=\sum_{\mu\in G} g_\mu$ with $f_\lambda\in E_\lambda$, $g_\mu\in E_\mu$ and
finite supports $F,G$.

Fix $\nu\in\Spec(\Delta)$. Since $U_\varphi$ acts on $E_\nu$ by multiplication by $e^{i\varphi(\nu)}$ and commutes with $P_\nu$,
\begin{align*}
P_\nu\!\left(U_\varphi^{-1}\big((U_\varphi f)(U_\varphi g)\big)\right)
&= e^{-i\varphi(\nu)}\,P_\nu\!\left(\Big(\sum_{\lambda\in F}e^{i\varphi(\lambda)}f_\lambda\Big)
\Big(\sum_{\mu\in G}e^{i\varphi(\mu)}g_\mu\Big)\right)\\
&= \sum_{\lambda\in F,\mu\in G} e^{i(\varphi(\lambda)+\varphi(\mu)-\varphi(\nu))}\,P_\nu(f_\lambda g_\mu)\\
&= \sum_{\lambda\in F,\mu\in G} \omega^{\nu}_{\lambda\mu}(\varphi)\,P_\nu(f_\lambda g_\mu)\\
&= P_\nu\big(f\star_{\omega(\varphi)} g\big),
\end{align*}
where the last line is exactly Definition~\ref{def:star} of $\star_{\omega(\varphi)}$ on $E_{\mathrm{fin}}$.
Since the $L^2$-spectral decomposition is orthogonal, equality of all spectral components implies \eqref{eq:gauge-conjugation} in $L^2(M)$.
Applying $U_\varphi$ to \eqref{eq:gauge-conjugation} gives
$U_\varphi(f\star_{\omega(\varphi)} g)=(U_\varphi f)(U_\varphi g)$ for all $f,g\in E_{\mathrm{fin}}$,
i.e.\ the claimed core conjugacy.
The algebra-isomorphism statement on $H^s$ then follows immediately from the fact that $U_\varphi$
is unitary on $H^s(M)$ and the product on $H^s$ is defined by the same conjugation formula.

Now fix $s>\frac{\dim M}{2}$. Pointwise multiplication is bounded $H^s\times H^s\to H^s$, and $U_\varphi$ commutes with $\Delta$,
hence is unitary on each Sobolev space $H^s(M)$. Therefore \eqref{eq:gauge-conjugation} defines a bounded bilinear map
$H^s\times H^s\to H^s$. The Fr\'echet algebra statement follows by restriction to $\bigcap_s H^s(M)=C^\infty(M)$.

The final $*$-compatibility sentence follows from Proposition~\ref{prop:star-id}.
\end{proof}

\begin{remark} Because pointwise multiplication is commutative, the gauge-twisted products $\star_{\omega(\varphi)}$ are commutative.
Producing intrinsic \emph{noncommutative} associative deformations without extra geometric input
appears to require richer channel spaces; this is part of the motivation for the future work outlined in \S\ref{sec:outlook}.
\end{remark}

\section{Sobolev regularity and associativity: general criteria}

The product $\star_\omega$ is always well-defined on $E_{\mathrm{fin}}$ with values in $L^2$ (Proposition~\ref{prop:L2-well}).
To iterate it and obtain a Banach $*$-algebra, one needs an ambient space closed under $\star_\omega$.
We isolate this analytically.

\begin{assumption}[Sobolev boundedness]\label{ass:A}
There exists $s>\frac{\dim M}{2}$ such that $\star_\omega$ extends to a bounded bilinear map
\[ \star_\omega:H^s(M)\times H^s(M)\to H^s(M). \]
\end{assumption}

\begin{lemma}[Uniqueness of bounded bilinear extension]\label{lem:unique-extension}
Let $X$ be a Banach space. If two bounded bilinear maps $B_1,B_2:H^s(M)\times H^s(M)\to X$ agree on
$E_{\mathrm{fin}}\times E_{\mathrm{fin}}$, then $B_1=B_2$.
\end{lemma}

\begin{proof}
Let $f,g\in H^s$. Choose sequences $f_n,g_n\in E_{\mathrm{fin}}$ converging to $f,g$ in $H^s$ (Lemma~\ref{lem:parseval-Hs}).
Bounded bilinearity implies $B_j(f,g)=\lim_n B_j(f_n,g_n)$, hence the limits agree.
\end{proof}

\begin{proposition}[Iterated products in $H^s$]\label{prop:iterated}
Assume Assumption~\ref{ass:A}. Then for every $n\ge 2$ and every $f_1,\dots,f_n\in H^s(M)$,
any iterated product formed using $\star_\omega$ is a well-defined element of $H^s(M)$.
Moreover each fixed parenthesization defines a continuous $n$-linear map $(H^s)^n\to H^s$.
\end{proposition}

\begin{proof}
Induct on $n$, using boundedness of $\star_\omega$.
\end{proof}

\begin{proposition}[$*$-structure on $H^s$]\label{prop:star-Hs}
Assume Assumption~\ref{ass:A} and that $\omega$ is $*$-compatible.
Then $\overline{f\star_\omega g}=\overline g\star_\omega\overline f$ holds for all $f,g\in H^s(M)$.
\end{proposition}

\begin{proof}
 The map $(f,g)\mapsto \overline{f\star_\omega g}-\overline g\star_\omega\overline f$ is continuous $H^s\times H^s\to H^s$ under Assumption~\ref{ass:A}.
It vanishes on the dense subset $E_{\mathrm{fin}}\times E_{\mathrm{fin}}$ by Proposition~\ref{prop:star-id}, hence vanishes identically.
\end{proof}

\subsection*{Associativity}

\begin{proposition}[Associativity reduces to a dense subspace]\label{prop:dense-assoc}
Assume Assumption~\ref{ass:A}. Then $\star_\omega$ is associative on $H^s(M)$ if and only if it is associative on $E_{\mathrm{fin}}$.
\end{proposition}

\begin{proof}
Define $A(f,g,h)=(f\star_\omega g)\star_\omega h - f\star_\omega(g\star_\omega h)$ for $f,g,h\in H^s$.
By Proposition~\ref{prop:iterated}, $A$ is a continuous trilinear map $H^s\times H^s\times H^s\to H^s$.
If $A$ vanishes on $E_{\mathrm{fin}}^3$, it vanishes on all of $(H^s)^3$ by density.
\end{proof}

\begin{lemma}[Triple-product component expansions]\label{lem:triple}
Assume Assumption~\ref{ass:A}. Let $f\in E_\lambda$, $g\in E_\mu$, $h\in E_\kappa$ and fix $\rho\in\Spec(\Delta)$.
Then in $E_\rho$ one has
\begin{align}\label{eq:left-triple}
P_\rho\big((f\star_\omega g)\star_\omega h\big)
&=\sum_{\nu\in\Spec(\Delta)} \omega^{\rho}_{\nu\kappa}\omega^{\nu}_{\lambda\mu}\, m^{\rho}_{\nu\kappa}\big(m^{\nu}_{\lambda\mu}(f,g),h\big),\\
\label{eq:right-triple}
P_\rho\big(f\star_\omega(g\star_\omega h)\big)
&=\sum_{\sigma\in\Spec(\Delta)} \omega^{\rho}_{\lambda\sigma}\omega^{\sigma}_{\mu\kappa}\, m^{\rho}_{\lambda\sigma}\big(f,m^{\sigma}_{\mu\kappa}(g,h)\big),
\end{align}
where both series converge in $E_\rho$.
\end{lemma}

\begin{proof}
We prove \eqref{eq:left-triple}; the other is analogous.
Since $f\star_\omega g\in H^s$, its spectral expansion converges unconditionally in $H^s$ (Lemma~\ref{lem:parseval-Hs}).
Let $S_N=\sum_{\nu\in F_N} P_\nu(f\star_\omega g)$ be finite partial sums over an increasing exhaustion $F_N$ of the spectrum.
Then $S_N\to f\star_\omega g$ in $H^s$. By continuity of $\star_\omega:H^s\times H^s\to H^s$,
$S_N\star_\omega h\to (f\star_\omega g)\star_\omega h$ in $H^s$, hence also $P_\rho(S_N\star_\omega h)\to P_\rho((f\star_\omega g)\star_\omega h)$ in $E_\rho$.
For each $N$,
\[ P_\rho(S_N\star_\omega h)=\sum_{\nu\in F_N} P_\rho(P_\nu(f\star_\omega g)\star_\omega h)=\sum_{\nu\in F_N} \omega^{\rho}_{\nu\kappa}\,m^{\rho}_{\nu\kappa}(P_\nu(f\star_\omega g),h). \]
Since $f\in E_\lambda$ and $g\in E_\mu$, Definition~\ref{def:star} gives $P_\nu(f\star_\omega g)=\omega^{\nu}_{\lambda\mu}m^{\nu}_{\lambda\mu}(f,g)$.
Substitute and pass to the limit $N\to\infty$.
\end{proof}

\begin{theorem}[Labelwise associativity criterion]\label{thm:assoc}
Assume Assumption~\ref{ass:A}.
Then $\star_\omega$ is associative on $H^s(M)$ if and only if for all
$\lambda,\mu,\kappa,\rho\in\Spec(\Delta)$ and all $f\in E_\lambda$, $g\in E_\mu$, $h\in E_\kappa$,
\begin{equation}\label{eq:assoc}
\sum_{\nu} \omega^{\rho}_{\nu\kappa}\omega^{\nu}_{\lambda\mu}\, m^{\rho}_{\nu\kappa}\big(m^{\nu}_{\lambda\mu}(f,g),h\big)
=\sum_{\sigma} \omega^{\rho}_{\lambda\sigma}\omega^{\sigma}_{\mu\kappa}\, m^{\rho}_{\lambda\sigma}\big(f,m^{\sigma}_{\mu\kappa}(g,h)\big)
\end{equation}
in $E_\rho$.
\end{theorem}

\begin{proof}
By Proposition~\ref{prop:dense-assoc}, associativity on $H^s$ is equivalent to associativity on $E_{\mathrm{fin}}$.
By trilinearity, it suffices to test associativity on eigenvectors.
Fix $\rho$.
Associativity implies
$P_\rho((f\star_\omega g)\star_\omega h)=P_\rho(f\star_\omega(g\star_\omega h))$.
Using Lemma~\ref{lem:triple} yields \eqref{eq:assoc}.
Conversely, if \eqref{eq:assoc} holds for all labels and eigenvectors, then
$P_\rho((f\star_\omega g)\star_\omega h-f\star_\omega(g\star_\omega h))=0$ for all $\rho$.
By Lemma~\ref{lem:parseval-Hs} this forces $(f\star_\omega g)\star_\omega h=f\star_\omega(g\star_\omega h)$ in $H^s$ on eigenvectors,
and hence on $E_{\mathrm{fin}}$ by trilinearity; then on $H^s$ by Proposition~\ref{prop:dense-assoc}.
\end{proof}
%%%%%%%%%%%%%%%%%%%%%%%%%%%%%%%%
%%%%%%%%%%%%%%%%%%%%%%%%%%%%%%%%%%%%%%%%%%%%%%%%%%%%%%%%%%%%%%%%%%%%%

\section{Compatibility with classical strict deformations}\label{sec:classical}

The framework above does not assume any symmetry of $M$.
This section is a \emph{compatibility check}: we recall three classical strict deformation procedures for actions of locally compact abelian groups
in their original formulations (Rieffel, Connes--Landi, Kasprzak), and then explain how, in situations where the group action has a discrete spectral decomposition---for instance, in the compact abelian case where the algebra decomposes into homogeneous subspaces indexed by characters of the acting group---their homogeneous-core products are recovered uniformly as refined instances of our spectral-channel twist. 
We do \emph{not} claim to recover these theories in full generality for arbitrary actions; rather, we identify the common graded mechanism
that interfaces with the refined spectral labels used in this paper.

\medskip
\noindent\textbf{A standing graded hypothesis for comparison on compact manifolds.}
Whenever a compact abelian group $K$ acts smoothly and isometrically on $(M,g)$, the induced action $\alpha:K\to\Aut(C^\infty(M))$
commutes with $\Delta$, hence each eigenspace $E_\lambda$ refines into $K$-isotypic subspaces.  Writing $\Gamma:=\widehat K$,
for $\gamma\in\Gamma$ set
\[
C^\infty(M)_\gamma:=\{f\in C^\infty(M):\alpha_t(f)=\gamma(t)f\ \forall t\in K\},
\qquad
E_{\lambda,\gamma}:=E_\lambda\cap C^\infty(M)_\gamma.
\]
Then $E_\lambda=\bigoplus_{\gamma\in\Gamma}E_{\lambda,\gamma}$ (finite sum for each fixed $\lambda$), and products satisfy
$C^\infty(M)_\gamma\cdot C^\infty(M)_\eta\subseteq C^\infty(M)_{\gamma\eta}$.

\medskip
\noindent\textbf{Refined labels.} For $\lambda\in\Spec(\Delta)$ and $\gamma\in\Gamma$ let $P_{\lambda,\gamma}:L^2(M)\to E_{\lambda,\gamma}$ be the $L^2$-orthogonal projection. Whenever we specify twisting coefficients $\omega$ on refined labels $(\lambda,\gamma)$, we apply Definition~\ref{def:star} verbatim with $P_\lambda$ replaced by $P_{\lambda,\gamma}$ (and with $E_{\mathrm{fin}}$ replaced by the algebraic direct sum $\bigoplus_{\lambda,\gamma}^{\mathrm{alg}}E_{\lambda,\gamma}$ of finite joint sums). We refer to the resulting product as ``Definition~\ref{def:star} with refined labels''.

\subsection{Rieffel deformation (actions of $\bbR^d$)}\label{subsec:Rieffel} 
Let $A$ be a (dense) $*$-algebra of smooth vectors for a strongly continuous action $\alpha:\bbR^d\to\Aut(A)$ by $*$-automorphisms
(e.g.\ $A=C^\infty(M)$ with a smooth action).  Fix a real skew-symmetric $d\times d$ matrix $J$.

Following \cite{Rieffel}, the (smooth) Rieffel-deformed product on $A$ is defined on smooth vectors by the oscillatory integral
\begin{equation}\label{eq:Rieffel-osc-fixed}
f\times_J g \ :=\ (2\pi)^{-d}\int_{\bbR^d}\int_{\bbR^d} \alpha_{Ju}(f)\,\alpha_v(g)\,
e^{-i\langle u,v\rangle}\,du\,dv,
\end{equation}
interpreted in the standard oscillatory sense (equivalently via regularization by Schwartz cutoffs).

To compare with our refined Laplace labels on a compact manifold, we impose the mild simplifying hypothesis that the $\bbR^d$-action is
\emph{periodic} and comes from an isometric action of a compact abelian group.
Concretely, assume $K=\bbT^d$ acts smoothly and isometrically on $M$ with induced action $\alpha^K$ on $C^\infty(M)$,
and let $\alpha$ be the lifted $\bbR^d$-action $\alpha_t:=\alpha^K_{t\!\!\mod 2\pi\bbZ^d}$.
Then $\Gamma=\widehat K\cong\bbZ^d$, and the algebraic direct sum $\bigoplus_{k\in\bbZ^d}^{\mathrm{alg}}C^\infty(M)_k$ of isotypic subspaces (finite Fourier sums) is dense in $C^\infty(M)$; for $f_k\in C^\infty(M)_k$ one has $\alpha_t(f_k)=e^{i\langle k,t\rangle}f_k$.

\begin{lemma}[Homogeneous multiplication]\label{lem:Rieffel-hom-fixed}
If $f_k\in C^\infty(M)_k$ and $g_\ell\in C^\infty(M)_\ell$ (with $k,\ell\in\bbZ^d$), then
\[
f_k\times_J g_\ell \ =\ \exp\!\big(i\langle k,J\ell\rangle\big)\, f_k g_\ell.
\]
\end{lemma}

\begin{proof}
Since $\alpha_{Ju}(f_k)=e^{i\langle k,Ju\rangle}f_k$ and $\alpha_v(g_\ell)=e^{i\langle \ell,v\rangle}g_\ell$,
the integrand in \eqref{eq:Rieffel-osc-fixed} equals
\[
e^{i\langle k,Ju\rangle}e^{i\langle \ell,v\rangle}e^{-i\langle u,v\rangle}\,f_k g_\ell
= e^{i\langle k,Ju\rangle}e^{i\langle \ell-u,v\rangle}\,f_k g_\ell.
\]
Fourier inversion gives $(2\pi)^{-d}\int_{\bbR^d} e^{i\langle \ell-u,v\rangle}\,dv=\delta(u-\ell)$ in the oscillatory sense,
so the $u$-integration collapses to $u=\ell$, yielding the stated phase factor.
\end{proof}

\paragraph{(D) Identification with our refined channel twist.}
Define a bicharacter $\sigma_J:\bbZ^d\times\bbZ^d\to\bbT$ by
\begin{equation}\label{eq:sigmaJ-fixed}
\sigma_J(k,\ell):=\exp\!\big(i\langle k,J\ell\rangle\big).
\end{equation}
Work with the refined decomposition $L^2(M)=\bigoplus_{\lambda,k}E_{\lambda,k}$ and define twisting coefficients on refined labels by
\[
\omega^{(\nu,k+\ell)}_{(\lambda,k)(\mu,\ell)} := \sigma_J(k,\ell).
\]

\begin{proposition}[Recovery of Rieffel deformation on the homogeneous core]\label{prop:recover-Rieffel-fixed}
Under the periodic/isometric hypothesis above, the product induced by Definition~\ref{def:star} (with refined labels) satisfies
\[
f_{\lambda,k}\star_\omega g_{\mu,\ell}=\sigma_J(k,\ell)\,f_{\lambda,k}g_{\mu,\ell}
\qquad (f_{\lambda,k}\in E_{\lambda,k},\ g_{\mu,\ell}\in E_{\mu,\ell}),
\]
and hence agrees on finite joint sums with the Rieffel product $f\times_J g$.
\end{proposition}

\begin{proof}
For homogeneous joint components $f_{\lambda,k}\in E_{\lambda,k}\subset C^\infty(M)_k$ and $g_{\mu,\ell}\in E_{\mu,\ell}\subset C^\infty(M)_\ell$,
Lemma~\ref{lem:Rieffel-hom-fixed} gives $f_{\lambda,k}\times_J g_{\mu,\ell}=\sigma_J(k,\ell)f_{\lambda,k}g_{\mu,\ell}$.
On the other hand, Definition~\ref{def:star} with the above $\omega$ yields
$P_{\nu,k+\ell}(f_{\lambda,k}\star_\omega g_{\mu,\ell})=\sigma_J(k,\ell)P_{\nu,k+\ell}(f_{\lambda,k}g_{\mu,\ell})$,
and summing over $\nu$ gives $f_{\lambda,k}\star_\omega g_{\mu,\ell}=\sigma_J(k,\ell)f_{\lambda,k}g_{\mu,\ell}$.
Extend by bilinearity to finite joint sums.
\end{proof}

\begin{remark}[Conventions]
Different normalizations of the oscillatory kernel in \eqref{eq:Rieffel-osc-fixed} lead to harmless rescalings of $J$ and/or factors of $2\pi$
in \eqref{eq:sigmaJ-fixed}. The identification above is at the level of the homogeneous-core phase factor.
\end{remark}

\subsection{Connes--Landi isospectral torus deformation}\label{subsec:CL}

Assume $\bbT^d$ acts smoothly and isometrically on $M$, with induced action $\alpha:\bbT^d\to\Aut(C^\infty(M))$ and isotypic subspaces $C^\infty(M)_k$.
The algebraic direct sum $\bigoplus_{k\in\bbZ^d}^{\mathrm{alg}}C^\infty(M)_k$ (finite Fourier sums) is dense in $C^\infty(M)$.

Fix a real skew-symmetric matrix $\Theta\in M_d(\bbR)$.
Let $C^\infty(\bbT^d_\Theta)$ denote the smooth noncommutative $d$-torus, i.e.\ the Fr\'echet $*$-algebra generated by unitaries
$U_1,\dots,U_d$ subject to the relations
\[
U_jU_k=e^{2\pi i\,\Theta_{jk}}U_kU_j,\qquad 1\le j,k\le d.
\]
For $k=(k_1,\dots,k_d)\in\bbZ^d$ define the (Weyl-normalized) Fourier unitaries
\[
U^k \;:=\; \exp\!\Bigl(-\pi i \sum_{1\le r<s\le d} k_r\,\Theta_{rs}\,k_s\Bigr)\,U_1^{k_1}\cdots U_d^{k_d}.
\]
Then a direct computation from the relations gives
\[
U^kU^\ell=\sigma_\Theta(k,\ell)\,U^{k+\ell},
\qquad 
\sigma_\Theta(k,\ell):=\exp(\pi i\,k^T\Theta\ell).
\]
Connes--Landi \cite{ConnesLandi} define the deformed smooth algebra as the fixed-point algebra
\[
C^\infty(M)_\Theta \ :=\ \bigl(C^\infty(M)\,\widehat\otimes\,C^\infty(\bbT^d_\Theta)\bigr)^{\bbT^d},
\]
for the diagonal action $t\cdot(f\otimes a):=\alpha_t(f)\otimes \gamma_{t^{-1}}(a)$, where $\gamma$ is the gauge action determined by $\gamma_t(U_j)=t_jU_j$ (equivalently, $\gamma_t(U^k)=t^kU^k$).

A dense $*$-subalgebra of the fixed-point algebra consists of finite Fourier sums of the form $\sum_k f_k\otimes U^k$
with $f_k\in C^\infty(M)_k$,
and multiplication in $C^\infty(M)_\Theta$ satisfies
\[
(f_k\otimes U^k)(g_\ell\otimes U^\ell)=\sigma_\Theta(k,\ell)\,f_kg_\ell\otimes U^{k+\ell}.
\]
Equivalently, transporting this product back to finite Fourier sums in $C^\infty(M)$ yields the homogeneous rule
\[
f_k\times_\Theta g_\ell \ =\ \sigma_\Theta(k,\ell)\,f_kg_\ell,
\]
and extension by bilinearity to finite sums.

Since the action is isometric, $E_\lambda=\bigoplus_k E_{\lambda,k}$ and we may work with refined labels $(\lambda,k)$.
Define twisting coefficients on refined labels by
\[
\omega^{(\nu,k+\ell)}_{(\lambda,k)(\mu,\ell)}:=\sigma_\Theta(k,\ell).
\]

\begin{proposition}[Recovery of the Connes--Landi product on the homogeneous core]\label{prop:recover-CL-fixed}
With $\omega$ as above, Definition~\ref{def:star} (with refined labels) satisfies
\[
f_{\lambda,k}\star_\omega g_{\mu,\ell}=\sigma_\Theta(k,\ell)\,f_{\lambda,k}g_{\mu,\ell}
\qquad (f_{\lambda,k}\in E_{\lambda,k},\ g_{\mu,\ell}\in E_{\mu,\ell}),
\]
so the resulting product agrees, on finite joint sums, with the transported Connes--Landi product $f\times_\Theta g$.
\end{proposition}

\begin{proof}
The proof is identical to the computation in Proposition~\ref{prop:recover-Rieffel-fixed} with $\sigma_J$ replaced by $\sigma_\Theta$:
for homogeneous joint components $f_{\lambda,k}$ and $g_{\mu,\ell}$ one has $f_{\lambda,k}g_{\mu,\ell}\in C^\infty(M)_{k+\ell}$,
and the refined spectral projections satisfy
$P_{\nu,k+\ell}(f_{\lambda,k}\star_\omega g_{\mu,\ell})=\sigma_\Theta(k,\ell)P_{\nu,k+\ell}(f_{\lambda,k}g_{\mu,\ell})$.
Summing over $\nu$ yields the stated homogeneous rule, and bilinearity gives the claim on finite sums.
\end{proof}

\subsection{Kasprzak deformation via Landstad theory (actions of locally compact abelian groups)}\label{subsec:Kasprzak}

Let $(A,\alpha)$ be a $C^*$-dynamical system for a locally compact abelian group $G$, so $\alpha:G\to\Aut(A)$ is strongly continuous,
and let $B:=A\rtimes_\alpha G$ be the crossed product with dual action $\widehat\alpha:\widehat G\to\Aut(B)$.

Let $\Omega:\widehat G\times\widehat G\to\bbT$ be a measurable normalized $2$-cocycle.
Kasprzak \cite{Kasprzak} constructs from $\Omega$ a cocycle-perturbed dual action $\widehat\alpha^\Omega$ on $B$,
and defines the deformed algebra $A^\Omega$ as the \emph{Landstad algebra} inside $M(B)$ associated to the $G$-product
$(B,\lambda,\widehat\alpha^\Omega)$. Equivalently, as the Landstad fixed-point subalgebra in $M(B)$ satisfying the usual Landstad conditions.
This yields a $C^*$-algebra $A^\Omega$ together with a canonical (deformed) $G$-action.

To interface with the refined Laplace labels on a compact manifold, we now specialize as follows:
assume $G=K$ is compact abelian and the action is induced by an isometric action on $(M,g)$, so that $\Gamma:=\widehat K$ is discrete and
\[
C^\infty(M)_{\mathrm{fin}}:=\bigoplus_{\gamma\in\Gamma}^{\mathrm{alg}}C^\infty(M)_\gamma\ \text{is dense in}\ C^\infty(M),\qquad E_\lambda=\bigoplus_{\gamma\in\Gamma}E_{\lambda,\gamma}.
\]
In this compact abelian setting, the cocycle deformation has a standard homogeneous description on the finite spectral core (see e.g.\ \cite{Kasprzak}):
there is a normalized unitary $2$-cocycle $\sigma:\Gamma\times\Gamma\to\bbT$ (determined by $\Omega$ together with the chosen cocycle-twist normalization)
such that, for homogeneous elements $a_\gamma\in C^\infty(M)_\gamma$ and $b_\eta\in C^\infty(M)_\eta$,
\begin{equation}\label{eq:Kasprzak-hom-fixed}
a_\gamma\times^\mathrm{Kas}_\Omega b_\eta \ =\ \sigma(\gamma,\eta)\,a_\gamma b_\eta.
\end{equation}

Work with the refined decomposition $L^2(M)=\bigoplus_{\lambda,\gamma}E_{\lambda,\gamma}$ and define twisting coefficients by
\[
\omega^{(\nu,\gamma\eta)}_{(\lambda,\gamma)(\mu,\eta)}:=\sigma(\gamma,\eta).
\]

\begin{proposition}[Recovery of the Kasprzak product on the homogeneous core]\label{prop:recover-Kasprzak-fixed}
In the compact abelian (discrete character-decomposition) setting above, Definition~\ref{def:star} (with refined labels) satisfies
\[
f_{\lambda,\gamma}\star_\omega g_{\mu,\eta}=\sigma(\gamma,\eta)\,f_{\lambda,\gamma}g_{\mu,\eta}
\qquad (f_{\lambda,\gamma}\in E_{\lambda,\gamma},\ g_{\mu,\eta}\in E_{\mu,\eta}),
\]
and hence agrees on finite joint sums with the Kasprzak product \eqref{eq:Kasprzak-hom-fixed} on smooth vectors.
\end{proposition}

\begin{proof}
The identity \eqref{eq:Kasprzak-hom-fixed} is the homogeneous component rule for the Kasprzak product in the compact abelian (discrete character-decomposition) case.
On the other hand, by construction of $\omega$ and the same refined-projection computation used above,
\[
P_{\nu,\gamma\eta}(f_{\lambda,\gamma}\star_\omega g_{\mu,\eta})
=\sigma(\gamma,\eta)\,P_{\nu,\gamma\eta}(f_{\lambda,\gamma}g_{\mu,\eta}),
\]
and summing over $\nu$ yields the homogeneous rule for $\star_\omega$.
Thus the two products agree on homogeneous elements and hence on finite joint sums by bilinearity.
\end{proof}

\begin{remark}
The comparison above isolates the common mechanism that appears when the relevant abelian group action has a discrete spectral decomposition, so that the algebra decomposes into homogeneous subspaces indexed by characters of the acting group.  This is the case for compact abelian actions and for periodic instances of $\bbR^d$-actions.  In this setting, our refined spectral-channel twisting reproduces the classical cocycle-deformed products uniformly on the corresponding homogeneous cores.
\end{remark}

%%%%%%%%%%%%%%%%%%%%%%%%%%%%%%%%%%%%%%%%%%%%%%%%%%%%%%%%%%%%%%%%%%%%%%%%%%%%%%%%%%%%%%%%%%%%%%%%%%%%%%
\section{New examples}

\subsection{An intrinsic action-free example: spectral gauge deformation}\label{subsec:intrinsic-gauge-example}

The deformations in \S\ref{sec:classical} require an additional abelian symmetry (a group action).
In contrast, the following family exists on \emph{every} compact Riemannian manifold and uses only the Laplace spectral decomposition.

Fix any function $\varphi:\Spec(\Delta)\to\bbR$ and define the spectral unitary
\[
U_\varphi\big|_{E_\lambda} = e^{i\varphi(\lambda)}\id_{E_\lambda}\qquad(\lambda\in\Spec(\Delta)).
\]
Equivalently, $U_\varphi = e^{i\varphi(\Delta)}$ by functional calculus.
Define a new product on $E_{\mathrm{fin}}$ by
\begin{equation}\label{eq:gauge-example-product}
f\star_\varphi g \ :=\ U_\varphi^{-1}\big((U_\varphi f)\,(U_\varphi g)\big),\qquad f,g\in E_{\mathrm{fin}}.
\end{equation}

This defines a bilinear map $E_{\mathrm{fin}}\times E_{\mathrm{fin}}\to L^2(M)$.
Moreover, by Theorem~\ref{thm:gauge-trivial}, for every $s>\frac{\dim M}{2}$ the same conjugation formula
extends to a strictly associative (indeed commutative) product on $H^s(M)$, and therefore restricts to the Fr\'echet algebra $C^\infty(M)=\bigcap_s H^s(M)$.

\begin{remark}
By Theorem~\ref{thm:gauge-trivial}, the product $\star_\varphi$ defined in \eqref{eq:gauge-example-product}
coincides with the scalar channel twist $\star_{\omega(\varphi)}$.
\end{remark}

\begin{example}[A concrete one-parameter family]\label{ex:gauge-parameter}
Taking $\varphi_t(\lambda):=t\log(1+\lambda)$ for $t\in\bbR$ gives $U_{\varphi_t}=(1+\Delta)^{it}$
and defines a one-parameter family of intrinsic products $\star_{\varphi_t}$.
This example uses no group actions and is not encompassed by the classical abelian-action deformation frameworks
(Rieffel/Connes--Landi/Kasprzak), though it is gauge-equivalent to the pointwise algebra by construction.
\end{example}

\begin{remark}
The intrinsic deformations above are ``spectral'' and action-free, but gauge-trivial (hence algebra-isomorphic to the undeformed algebra; and $\ast$-isomorphic either after transporting the involution or in the real-phase case).
Producing intrinsic \emph{noncommutative} strictly associative examples appears to require richer (higher-dimensional) channel spaces
and non-scalar recoupling/associator data, which is precisely the focus of the planned sequel.
\end{remark}

\section{Obstruction results}\label{sec:obstruction}

The identity in Theorem~\ref{thm:gauge-trivial} is a conjugacy of bilinear maps on $E_{\mathrm{fin}}$
(with values in $L^2(M)$), and for every $s>\frac{\dim M}{2}$ it upgrades to an algebra conjugacy on the Sobolev algebra $H^s(M)$. Since the same conjugacy holds for every $s$,
it restricts to the Fr\'echet algebra $C^\infty(M)=\bigcap_s H^s(M)$. If one equips $H^s(M)$ (for a fixed $s>\frac{\dim M}{2}$) with the transported involution 
\begin{equation}\label{eq:star_transport}
f^{\sharp}:=U_\varphi^{-1}\big(\overline{U_\varphi f}\big)    
\end{equation}
then $U_\varphi$ becomes a $*$-isomorphism onto $(H^s(M),\cdot, \overline{\phantom{f}})$,
and the same holds on $E_{\mathrm{fin}}$ by restriction.

If instead one insists on the \emph{usual} involution $f\mapsto\overline f$ on the deformed algebra,
then $\star_{\omega(\varphi)}$ is a $*$-algebra product only in the special case when the phases are real,
equivalently $e^{i\varphi(\lambda)}\in\{\pm1\}$ for all $\lambda$ so that $U_\varphi$ commutes with complex conjugation.

\begin{proposition}[Gauge equivalence by a spectral coboundary]\label{prop:gauge_equivalence}
Let $\omega=(\omega^\nu_{\lambda\mu})$ and $\omega'=(\omega'{}^\nu_{\lambda\mu})$ be two scalar twisting systems.
Assume there exists a function $a:\Spec(\Delta)\to\bbT$ such that for all $\lambda,\mu,\nu$,
\begin{equation}\label{eq:coboundary_relation}
\omega'{}^\nu_{\lambda\mu} \;=\; a(\lambda)a(\mu)\overline{a(\nu)}\,\omega^\nu_{\lambda\mu}.
\end{equation}
Define the diagonal unitary $A:=\sum_\lambda a(\lambda)P_\lambda$ on $L^2(M)$.

Then for all $f,g\in E_{\mathrm{fin}}$,
\[
A(f\star_{\omega'} g) \;=\; (Af)\star_\omega (Ag).
\]
Thus $A$ is a linear bijection $E_{\mathrm{fin}}\to E_{\mathrm{fin}}$ intertwining the core products (as maps into $L^2(M)$).
Moreover, if $\star_\omega$ and $\star_{\omega'}$ extend continuously to $H^s(M)$ for some $s$,
then $A$ extends to a bounded algebra isomorphism on $H^s(M)$ intertwining the extended products.
\end{proposition}

\begin{proof}
Write $f=\sum_\lambda f_\lambda$ and $g=\sum_\mu g_\mu$ with $f_\lambda\in E_\lambda$, $g_\mu\in E_\mu$ and finite supports.
For fixed $\nu$,
\[
P_\nu(f\star_{\omega'} g)=\sum_{\lambda,\mu} \omega'{}^\nu_{\lambda\mu}\,P_\nu(f_\lambda g_\mu).
\]
Applying $A$ multiplies $E_\nu$ by $a(\nu)$, so
\[
P_\nu\big(A(f\star_{\omega'} g)\big)
=
a(\nu)\sum_{\lambda,\mu} \omega'{}^\nu_{\lambda\mu}\,P_\nu(f_\lambda g_\mu).
\]
Using \eqref{eq:coboundary_relation} and $\overline{a(\nu)}=a(\nu)^{-1}$ gives
\[
P_\nu\big(A(f\star_{\omega'} g)\big)
=
\sum_{\lambda,\mu} a(\lambda)a(\mu)\,\omega^\nu_{\lambda\mu}\,P_\nu(f_\lambda g_\mu).
\]
On the other hand, $Af=\sum_\lambda a(\lambda)f_\lambda$ and $Ag=\sum_\mu a(\mu)g_\mu$, hence
\[
P_\nu\big((Af)\star_\omega (Ag)\big)
=
\sum_{\lambda,\mu} \omega^\nu_{\lambda\mu}\,P_\nu\big(a(\lambda)f_\lambda\cdot a(\mu)g_\mu\big)
=
\sum_{\lambda,\mu} a(\lambda)a(\mu)\,\omega^\nu_{\lambda\mu}\,P_\nu(f_\lambda g_\mu),
\]
which matches $P_\nu(A(f\star_{\omega'} g))$ for each $\nu$. Summing over $\nu$ in $L^2$ yields the claimed identity.

Finally, $A$ is unitary on $L^2$ and diagonal in the Laplace decomposition, hence bounded on every $H^s(M)$.
Therefore, if $\star_\omega$ and $\star_{\omega'}$ extend to $H^s$, the intertwining identity extends by continuity.
\end{proof}

\begin{remark}[$*$-isomorphisms]\label{rem:gauge_star_iso}
The map $A$ intertwines the core products on $E_{\mathrm{fin}}$ (as maps into $L^2(M)$).
Whenever the products extend to some $H^s(M)$, $A$ is an algebra isomorphism on $H^s(M)$.
It is a $*$-isomorphism for the standard involution $f\mapsto\overline f$ precisely when $A$ commutes with complex conjugation,
i.e.\ $a(\lambda)\in\{\pm1\}$ for all $\lambda$.
Otherwise, $A$ still becomes a $*$-isomorphism if one transports the involution as in ~(\ref{eq:star_transport}).
\end{remark}

\begin{proposition}[Rigidity of separable twists]\label{prop:separable_twists_are_coboundaries}
Assume $\omega=(\omega^\nu_{\lambda\mu})$ has the separable form
\[
\omega^\nu_{\lambda\mu}=p(\lambda)\,q(\mu)\,\overline{r(\nu)},
\qquad p,q,r:\Spec(\Delta)\to\bbT.
\]
Assume moreover that $\star_\omega$ is unital with unit $1\in E_0$, i.e.\ $f\star_\omega 1=f=1\star_\omega f$ for all $f\in E_{\mathrm{fin}}$.
Then $\omega$ is a spectral coboundary: there exists $a:\Spec(\Delta)\to\bbT$ with $a(0)=1$ such that
\[
\omega^\nu_{\lambda\mu}=a(\lambda)a(\mu)\overline{a(\nu)}.
\]
Consequently, $\star_\omega$ is gauge-trivial in the sense of Proposition~\ref{prop:gauge_equivalence}.
\end{proposition}

\begin{proof}
Unitality implies, for each $\lambda$, that $f_\lambda\star_\omega 1=f_\lambda$ for all $f_\lambda\in E_\lambda$.
Since $1\in E_0$ and $f_\lambda\cdot 1=f_\lambda\in E_\lambda$, we have
\[
P_\nu(f_\lambda\cdot 1)=\delta_{\nu\lambda}\,f_\lambda,
\]
so in the defining formula for $P_\nu(f_\lambda\star_\omega 1)$ only the term with $\nu=\lambda$ survives, and hence $\omega^\lambda_{\lambda 0}=1$.
By the separable form this reads
\[
p(\lambda)\,q(0)\,\overline{r(\lambda)}=1\qquad\text{for all }\lambda,
\]
hence $p(\lambda)=\overline{q(0)}\,r(\lambda)$.

Similarly, $1\star_\omega g_\mu=g_\mu$ for $g_\mu\in E_\mu$ gives $\omega^\mu_{0\mu}=1$, i.e.
\[
p(0)\,q(\mu)\,\overline{r(\mu)}=1\qquad\text{for all }\mu,
\]
so $q(\mu)=\overline{p(0)}\,r(\mu)$.

Therefore
\[
\omega^\nu_{\lambda\mu}
=
p(\lambda)\,q(\mu)\,\overline{r(\nu)}
=
\overline{q(0)p(0)}\ r(\lambda)r(\mu)\overline{r(\nu)}.
\]
Now $1\star_\omega 1=1$ implies $\omega^0_{00}=1$, so
\[
1=\omega^0_{00}=\overline{q(0)p(0)}\ r(0),
\qquad\text{hence }\overline{q(0)p(0)}=\overline{r(0)}.
\]
Thus
\[
\omega^\nu_{\lambda\mu}=\overline{r(0)}\ r(\lambda)r(\mu)\overline{r(\nu)}.
\]
Define $a:\Spec(\Delta)\to\bbT$ by $a(\lambda):=\overline{r(0)}\,r(\lambda)$.
Then $a(0)=\overline{r(0)}r(0)=1$ and
\[
a(\lambda)a(\mu)\overline{a(\nu)}
=
(\overline{r(0)}r(\lambda))(\overline{r(0)}r(\mu))\overline{\overline{r(0)}r(\nu)}
=
\overline{r(0)}\ r(\lambda)r(\mu)\overline{r(\nu)}
=
\omega^\nu_{\lambda\mu},
\]
as required.
\end{proof}

\subsection{A cohomological obstruction for graded scalar twists}\label{subsec:cohomological-obstruction}

In the scalar regime one can obtain a large, rigid, and fully classifiable class of strictly associative twists
when the twisting data is controlled by an additive grading compatible with the multiplication channels.
We first formulate this on a dense graded core determined by a refinement, where the cohomological classification is purely algebraic.
We then record the additional explicit continuity requirements needed to transport the classification to Sobolev completions.

\subsubsection*{Refined pieces and channel-compatible gradings on a Sobolev algebra}

Fix $s>\dim(M)/2$ so that $H^s(M)$ is a Banach algebra under pointwise multiplication.
Assume we are given a countable family of finite-dimensional subspaces $(V_\alpha)_{\alpha\in I}\subset H^s(M)$
such that
\begin{equation}\label{eq:Hs-refinement}
V_\alpha\perp V_\beta\ \text{in }L^2(M)\ \text{for }\alpha\neq\beta.
\end{equation}
Let $P_\alpha:L^2(M)\to V_\alpha$ denote the $L^2$-orthogonal projection.

Define refined multiplication channels
\[
m^\gamma_{\alpha\beta}:V_\alpha\times V_\beta\to V_\gamma,\qquad
m^\gamma_{\alpha\beta}(f,g):=P_\gamma(fg).
\]
We write $\alpha\otimes\beta\rightsquigarrow\gamma$ if $m^\gamma_{\alpha\beta}$ is not identically zero.

\begin{definition}[Universal grading group]\label{def:Gamma-univ}
Let $F_I:=\bigoplus_{\alpha\in I}\mathbb{Z}e_\alpha$ be the free abelian group on generators $e_\alpha$.
Let $R\subseteq F_I$ be the subgroup generated by all $e_\alpha+e_\beta-e_\gamma$ such that
$\alpha\otimes\beta\rightsquigarrow\gamma$.
The \emph{universal grading group} of the refinement $(V_\alpha)$ is
\[
\Gamma_{\mathrm{univ}}:=F_I/R,
\]
and we write $[\alpha]\in\Gamma_{\mathrm{univ}}$ for the class of $e_\alpha$.
\end{definition}

A map $d:I\to\Gamma$ into an abelian group is called \emph{channel-compatible} if
\[
d(\gamma)=d(\alpha)+d(\beta)\qquad\text{whenever}\qquad \alpha\otimes\beta\rightsquigarrow\gamma.
\]
Equivalently, $d$ is the composite of the canonical map $\alpha\mapsto[\alpha]\in\Gamma_{\mathrm{univ}}$
with a unique homomorphism $\Gamma_{\mathrm{univ}}\to\Gamma$.

\begin{definition}[Degree subspaces and finite-degree Sobolev core]\label{def:Hs-degree-core}
Assume $s>\dim(M)/2$ so that $H^s(M)$ is a Banach algebra under pointwise multiplication.
Let $(V_\alpha)_{\alpha\in I}$ be an $L^2$-orthogonal refinement of the Laplace eigenspaces with
$V_\alpha\subseteq E_{\lambda(\alpha)}\subset C^\infty(M)$ and $\bigoplus_{\alpha}^{\mathrm{alg}}V_\alpha=E_{\mathrm{fin}}$.
Let $d:I\to\Gamma$ be a channel-compatible map to an abelian group $\Gamma$.

For $p\in\Gamma$ define the closed $L^2$-subspace
\[
L^2_p:=\overline{\bigoplus_{\alpha:\ d(\alpha)=p}^{\mathrm{alg}}V_\alpha}^{\,L^2},
\]
and define the degree-$p$ Sobolev subspace
\[
H^s_p:=H^s(M)\cap L^2_p,
\]
equipped with the subspace norm from $H^s(M)$.

Finally, define the finite-degree Sobolev core
\[
H^s_{\Gamma,\mathrm{fin}}:=\bigoplus_{p\in\Gamma}^{\mathrm{alg}}H^s_p,
\]
i.e.\ the subspace of elements having finite support in $\Gamma$. We write $(H^s_{\Gamma,\mathrm{fin}})_p:=H^s_p\subset H^s_{\Gamma,\mathrm{fin}}$ for the degree-$p$ summand.
\end{definition}

\begin{lemma}[Multiplication preserves degree]\label{lem:Hs-degree-multiplication}
Under the hypotheses of Definition~\ref{def:Hs-degree-core}, for all $p,q\in\Gamma$ one has
\[
H^s_p\cdot H^s_q \subseteq H^s_{p+q}.
\]
Consequently, $H^s_{\Gamma,\mathrm{fin}}$ is a graded subalgebra of $H^s(M)$, and it is dense in $H^s(M)$.
\end{lemma}

\begin{proof}
Let $p,q\in\Gamma$.
First consider $f\in \bigoplus_{d(\alpha)=p}^{\mathrm{alg}}V_\alpha$ and
$g\in \bigoplus_{d(\beta)=q}^{\mathrm{alg}}V_\beta$.
By channel-compatibility of $d$, every nonzero $V_\gamma$-component of $fg$ must satisfy $d(\gamma)=p+q$.
Since the $V_\gamma$ form an $L^2$-orthogonal refinement, this implies $fg\in L^2_{p+q}$.
Also $f,g\in H^s(M)$ and $H^s(M)$ is an algebra, hence $fg\in H^s(M)$, so $fg\in H^s_{p+q}$.

Now let $f\in H^s_p$ and $g\in H^s_q$.
Let $\Pi_N:=\sum_{\lambda\le N}P_\lambda$ be the Laplace spectral cutoff.
Then $\Pi_N$ is bounded on $H^s(M)$ and $\Pi_N\to \mathrm{id}$ strongly on $H^s(M)$.
Moreover, $\Pi_N$ preserves degrees because it is diagonal with respect to the refinement $(V_\alpha)$, hence
$\Pi_N f \in \bigoplus_{d(\alpha)=p}^{\mathrm{alg}}V_\alpha$ and
$\Pi_N g \in \bigoplus_{d(\beta)=q}^{\mathrm{alg}}V_\beta$.
By the previous paragraph, $(\Pi_N f)(\Pi_N g)\in H^s_{p+q}$ for each $N$.
Since multiplication is continuous on $H^s(M)$, $(\Pi_N f)(\Pi_N g)\to fg$ in $H^s(M)$.
Because $H^s_{p+q}$ is closed in $H^s(M)$ (being $H^s(M)\cap L^2_{p+q}$ with $L^2_{p+q}$ closed),
we conclude $fg\in H^s_{p+q}$.

Finally, $H^s_{\Gamma,\mathrm{fin}}$ is closed under multiplication by the degree rule and finite support,
and it is dense since it contains $E_{\mathrm{fin}}=\bigoplus_\alpha^{\mathrm{alg}}V_\alpha$, which is dense in $H^s(M)$.
\end{proof}

\subsubsection*{Cocycle twists and the obstruction class}

\begin{definition}[Normalized unitary $2$-cocycle]\label{def:2cocycle}
A map $\sigma:\Gamma\times\Gamma\to\mathbb{T}$ is a normalized unitary $2$-cocycle if
\[
\sigma(p,0)=\sigma(0,p)=1
\quad\text{and}\quad
\sigma(p,q)\,\sigma(p+q,r)=\sigma(q,r)\,\sigma(p,q+r)
\]
for all $p,q,r\in\Gamma$.
Two cocycles $\sigma,\sigma'$ are \emph{cohomologous} if there exists $b:\Gamma\to\mathbb{T}$ with
$\sigma'(p,q)=b(p)b(q)\overline{b(p+q)}\,\sigma(p,q)$.
\end{definition}

\begin{definition}[Graded cocycle twisted product on the dense graded core]\label{def:graded-cocycle-product}
Let $d:I\to\Gamma$ be channel-compatible and let $\sigma$ be a normalized unitary $2$-cocycle.
For homogeneous $f\in (H^s_{\Gamma,\mathrm{fin}})_p$ and $g\in (H^s_{\Gamma,\mathrm{fin}})_q$ define
\[
f\star_{\sigma} g:=\sigma(p,q)\,fg,
\]
and extend by bilinearity to $H^s_{\Gamma,\mathrm{fin}}\times H^s_{\Gamma,\mathrm{fin}}$.

Whenever $\star_\sigma$ extends (uniquely) to a bounded bilinear map $H^s(M)\times H^s(M)\to H^s(M)$, we use the same notation
for the extension.
\end{definition}

\begin{proposition}[Associativity and commutativity on the graded core]\label{prop:graded-assoc}
The product $\star_\sigma$ on $H^s_{\Gamma,\mathrm{fin}}$ is strictly associative.
Moreover, $\star_\sigma$ is commutative on $H^s_{\Gamma,\mathrm{fin}}$ if and only if the commutator bicharacter
\[
\beta_\sigma(p,q):=\sigma(p,q)\,\overline{\sigma(q,p)}
\]
is identically $1$ on the subgroup of $\Gamma$ generated by $\{d(\alpha):\alpha\in I\}$.

If $\star_\sigma$ extends to a bounded bilinear map $H^s\times H^s\to H^s$, then the extension is associative (and commutative
under the same condition on $\beta_\sigma$).
\end{proposition}

\begin{proof}
We first verify associativity on the algebraic graded core $H^s_{\Gamma,\mathrm{fin}}$.
Let $f\in (H^s_{\Gamma,\mathrm{fin}})_p$, $g\in (H^s_{\Gamma,\mathrm{fin}})_q$, and $h\in (H^s_{\Gamma,\mathrm{fin}})_r$ be homogeneous.
Then, by Definition~\ref{def:graded-cocycle-product},
\[
(f\star_\sigma g)\star_\sigma h
=\sigma(p,q)\sigma(p+q,r)\,fgh,
\qquad
f\star_\sigma(g\star_\sigma h)
=\sigma(q,r)\sigma(p,q+r)\,fgh,
\]
and these coincide by the cocycle identity in Definition~\ref{def:2cocycle}. Hence $\star_\sigma$ is associative on homogeneous elements,
and therefore on $H^s_{\Gamma,\mathrm{fin}}$ by bilinearity.

For commutativity on homogeneous elements one has
\[
f\star_\sigma g=\sigma(p,q)\,fg,\qquad
g\star_\sigma f=\sigma(q,p)\,gf=\sigma(q,p)\,fg,
\]
so $f\star_\sigma g=g\star_\sigma f$ for all homogeneous $f,g$ if and only if $\beta_\sigma(p,q)=1$
for all degrees $p,q$ that occur.  Since (by a routine cocycle computation) $\beta_\sigma$ is a bicharacter on $\Gamma$,
this is equivalent to $\beta_\sigma\equiv 1$ on the subgroup of $\Gamma$ generated by $\{d(\alpha):\alpha\in I\}$.

Finally, suppose $\star_\sigma$ extends to a bounded bilinear map $H^s\times H^s\to H^s$.
Then the maps
\[
(f,g,h)\longmapsto (f\star_\sigma g)\star_\sigma h,
\qquad
(f,g,h)\longmapsto f\star_\sigma(g\star_\sigma h)
\]
are bounded trilinear maps on $(H^s)^3$.
They agree on the dense subset $(H^s_{\Gamma,\mathrm{fin}})^3$ by the associativity already proved on $H^s_{\Gamma,\mathrm{fin}}$,
hence they agree on all of $(H^s)^3$ by continuity. The commutativity statement upgrades in the same way.
\end{proof}

\begin{corollary}[Rank-one obstruction to noncommutativity]\label{cor:rank-one}
Let $\Gamma_0\le \Gamma$ be the subgroup generated by the degrees that occur, i.e.\ $\Gamma_0:=\langle d(\alpha):\alpha\in I\rangle$.
If $\Gamma_0$ is cyclic (e.g.\ $\Gamma_0\cong\bbZ$ or a finite cyclic group), then for every normalized unitary $2$-cocycle $\sigma$ on $\Gamma$ the commutator bicharacter
$\beta_\sigma$ is trivial on $\Gamma_0$, and hence the graded cocycle twist $\star_\sigma$ is commutative on $H^s_{\Gamma,\mathrm{fin}}$
(and on any Sobolev completion where $\star_\sigma$ extends).
\end{corollary}

\begin{proof}
By Proposition~\ref{prop:graded-assoc}, $\beta_\sigma$ is a bicharacter on $\Gamma$ and satisfies $\beta_\sigma(p,p)=1$ for all $p\in\Gamma$.
If $\Gamma_0=\langle p_0\rangle$ is cyclic, then for all $m,n\in\bbZ$,
\[
\beta_\sigma(mp_0,np_0)=\beta_\sigma(p_0,p_0)^{mn}=1,
\]
so $\beta_\sigma\equiv 1$ on $\Gamma_0$. The commutativity conclusion now follows again from Proposition~\ref{prop:graded-assoc}.
\end{proof}

For a $1$-cochain $b:\Gamma\to\bbT$, we write
\[
U_b\big|_{(H^s_{\Gamma,\mathrm{fin}})_p} := b(p)\,\mathrm{id},
\qquad p\in\Gamma,
\]
for the diagonal unitary on $H^s_{\Gamma,\mathrm{fin}}$ acting by phase $b(p)$ on the degree-$p$ component, that is, $U_b(f_p)=b(p)\,f_p$ for $f_p\in H^s_p$.

Recall that the second cohomology group $H^2(\Gamma,\bbT)$ is defined \cite{BrownCohomology} as the quotient
\[
H^2(\Gamma,\bbT)
:=
Z^2(\Gamma,\bbT)\big/ B^2(\Gamma,\bbT),
\]
where $Z^2(\Gamma,\bbT)$ consists of \emph{normalized} $2$-cocycles $\sigma:\Gamma\times\Gamma\to\bbT$
satisfying
\[
\sigma(p,0)=\sigma(0,p)=1
\quad\text{and}\quad
\sigma(p,q)\sigma(p+q,r)=\sigma(q,r)\sigma(p,q+r),
\]
and $B^2(\Gamma,\bbT)$ consists of $2$-coboundaries of the form
\[
\sigma(p,q)=b(p)b(q)\overline{b(p+q)}
\]
for $1$-cochains $b:\Gamma\to\bbT$.

\begin{theorem}[Cohomological classification and obstruction]\label{thm:cohomological-obstruction}
Fix a refinement $(V_\alpha)$ as in \eqref{eq:Hs-refinement} and a channel-compatible grading $d:I\to\Gamma$.

\smallskip
\noindent (i) On the dense graded core $H^s_{\Gamma,\mathrm{fin}}$, graded cocycle twisted products $\star_\sigma$ are classified up to diagonal graded gauge
by the cohomology group $H^2(\Gamma,\bbT)$: if $\sigma'$ is cohomologous to $\sigma$ via a $1$-cochain $b:\Gamma\to\bbT$, i.e.
\[
\sigma'(p,q)=b(p)b(q)\overline{b(p+q)}\,\sigma(p,q),
\]
then
\[
f\star_{\sigma'} g \ =\ U_b^{-1}\bigl((U_b f)\star_\sigma (U_b g)\bigr)\qquad (f,g\in H^s_{\Gamma,\mathrm{fin}}),
\]
and conversely any such diagonal graded conjugacy arises from a coboundary.

\smallskip
\noindent (ii) (\emph{Sobolev upgrade.})
If, in addition, for a given pair of cocycles $\sigma,\sigma'$ the products $\star_\sigma,\star_{\sigma'}$ extend to bounded bilinear maps
$H^s\times H^s\to H^s$ \emph{and} the diagonal map $U_b$ extends to a bounded invertible operator on $H^s$,
then the same conjugacy identity holds on $H^s$ and the corresponding Sobolev deformations are isomorphic.

\smallskip
\noindent In particular, if $H^2(\Gamma,\bbT)=0$, then every graded cocycle twist in this channel-compatible subclass is diagonally gauge-trivial
on the graded core (and likewise on any completion where the above extension hypotheses hold).  Thus, within this subclass, there is no strictly associative
noncommutative deformation whenever $\beta_\sigma\equiv 1$ for every cocycle class (e.g.\ when $H^2(\Gamma,\bbT)=0$).
\end{theorem}

\begin{proof}
(i) If $\sigma'$ is cohomologous to $\sigma$ via $b$, define $U_b$ on $H^s_{\Gamma,\mathrm{fin}}$ by $U_b|_{(H^s_{\Gamma,\mathrm{fin}})_p}=b(p)\mathrm{id}$.
For homogeneous $f\in(H^s_{\Gamma,\mathrm{fin}})_p$ and $g\in(H^s_{\Gamma,\mathrm{fin}})_q$,
\[
U_b^{-1}\bigl((U_b f)\star_\sigma (U_b g)\bigr)
=\overline{b(p+q)}\,\sigma(p,q)\,b(p)b(q)\,fg
=\sigma'(p,q)\,fg
=f\star_{\sigma'}g,
\]
and bilinearity gives the identity on $H^s_{\Gamma,\mathrm{fin}}$.
Conversely, any diagonal graded conjugacy must be given by scalars on each homogeneous component, hence determines a $1$-cochain $b$
and forces the cocycles to differ by the corresponding coboundary.

(ii) If the stated extensions exist on $H^s$, then the conjugacy identity continues to hold by density of $H^s_{\Gamma,\mathrm{fin}}$ in $H^s$ and continuity.
The final statements are immediate consequences.
\end{proof}
\begin{remark}[Intrinsic meaning]\label{rem:intrinsic-meaning}
The group $\Gamma_{\mathrm{univ}}$ depends only on the intrinsic multiplication--projection pattern of the refinement,
hence is an isometry-invariant of $(M,g)$ once a canonical refinement has been specified (e.g.\ by joint eigenspaces in
the presence of additional commuting operators).  Theorem~\ref{thm:cohomological-obstruction} should be read as an
obstruction theorem for a large, structurally rigid subclass of scalar twists: those controlled by a channel-compatible
grading.
\end{remark}

%%%%%%%%%%%%%%%%%%%%%%%%%%%%%%%%%%%%%%%%%%%%%%%%%%%%%%%%%%%%%%%%%%%%%%
\section{From actions to gradings to cocycle-twisted products}\label{subsec:action-grading-twist}

The following ``grading principle'' makes precise the slogan that, in cocycle twisting,
the grading is the algebraic essence and a compact abelian action is a mechanism for producing that grading.

\begin{definition}[$\Gamma$-graded $*$-algebra]\label{def:graded-star-alg}
Let $\Gamma$ be a discrete abelian group (written additively).
A \emph{$\Gamma$-graded $*$-algebra} is a $*$-algebra $A$ together with a vector space decomposition
\[
A=\bigoplus_{\gamma\in\Gamma}^{\mathrm{alg}} A_\gamma
\]
such that $A_\gamma A_\eta\subseteq A_{\gamma+\eta}$ and $A_\gamma^*\subseteq A_{-\gamma}$ for all $\gamma,\eta\in\Gamma$.
\end{definition}

\begin{theorem}[Grading principle: actions supply gradings, cocycle twisting uses only the grading]\label{thm:A}
Let $K$ be a compact abelian group with discrete dual group $\Gamma=\widehat K$ (written additively), and let
$\alpha:K\to\mathrm{Aut}(A)$ be an action by $*$-automorphisms on a $*$-algebra $A$.
For each $\gamma\in\Gamma$ define the spectral subspace
\[
A_\gamma^\alpha:=\{a\in A:\ \alpha_t(a)=\gamma(t)\,a\ \ \forall t\in K\},
\]
and let
\[
A_{\mathrm{fin}}^\alpha:=\bigoplus_{\gamma\in\Gamma}^{\mathrm{alg}}A_\gamma^\alpha
\]
denote the finite spectral core.

\smallskip
\noindent (i) The decomposition $A_{\mathrm{fin}}^\alpha=\bigoplus_{\gamma}^{\mathrm{alg}}A_\gamma^\alpha$ is a $\Gamma$-grading:
\[
A_\gamma^\alpha A_\eta^\alpha\subseteq A_{\gamma+\eta}^\alpha,\qquad (A_\gamma^\alpha)^*\subseteq A_{-\gamma}^\alpha.
\]

\smallskip
\noindent (ii) Given a normalized unitary $2$-cocycle $\sigma$ on $\Gamma$ (Definition~\ref{def:2cocycle}), define a bilinear product on
$A_{\mathrm{fin}}^\alpha$ by the homogeneous rule
\begin{equation}\label{eq:graded-product}
a_\gamma\star_\sigma b_\eta := \sigma(\gamma,\eta)\,a_\gamma b_\eta,
\qquad a_\gamma\in A_\gamma^\alpha,\ b_\eta\in A_\eta^\alpha,
\end{equation}
extended by bilinearity to finite sums. Then $\star_\sigma$ is associative on $A_{\mathrm{fin}}^\alpha$.
Moreover, with the \emph{original} involution on $A$, the identity
\[
(a\star_\sigma b)^*=b^*\star_\sigma a^*\qquad(a,b\in A_{\mathrm{fin}}^\alpha)
\]
holds if and only if
\begin{equation}\label{eq:cocycle-star-compat}
\overline{\sigma(\gamma,\eta)}=\sigma(-\eta,-\gamma)\qquad (\gamma,\eta\in\Gamma).
\end{equation}

\smallskip
\noindent (iii) The product \eqref{eq:graded-product} depends only on the graded pieces $(A_\gamma^\alpha)_\gamma$ (and $\sigma$):
if $\beta:K\to\mathrm{Aut}(A)$ is another action with $A_\gamma^\beta=A_\gamma^\alpha$ for all $\gamma$ (same labeling),
then the corresponding cocycle-twisted products coincide on the common graded core.
\end{theorem}

\begin{proof}
(i) Each $A_\gamma^\alpha$ is a linear subspace, since each $\alpha_t$ is linear.
If $a\in A_\gamma^\alpha$ and $b\in A_\eta^\alpha$, then using $(\gamma+\eta)(t)=\gamma(t)\eta(t)$ one has
\[
\alpha_t(ab)=\alpha_t(a)\alpha_t(b)=\gamma(t)a\cdot\eta(t)b=(\gamma+\eta)(t)\,ab,
\]
so $ab\in A_{\gamma+\eta}^\alpha$. Also,
\[
\alpha_t(a^*)=\alpha_t(a)^*=\overline{\gamma(t)}\,a^*=(-\gamma)(t)\,a^*,
\]
so $a^*\in A_{-\gamma}^\alpha$.

(ii) Associativity reduces to homogeneous elements. For $a_\gamma\in A_\gamma^\alpha$, $b_\eta\in A_\eta^\alpha$, $c_\rho\in A_\rho^\alpha$,
\[
(a_\gamma\star_\sigma b_\eta)\star_\sigma c_\rho
=\sigma(\gamma,\eta)\sigma(\gamma+\eta,\rho)\,a_\gamma b_\eta c_\rho,
\qquad
a_\gamma\star_\sigma(b_\eta\star_\sigma c_\rho)
=\sigma(\eta,\rho)\sigma(\gamma,\eta+\rho)\,a_\gamma b_\eta c_\rho,
\]
and these coincide by the cocycle identity in Definition~\ref{def:2cocycle}.
For the involution, for homogeneous $a_\gamma,b_\eta$,
\[
(a_\gamma\star_\sigma b_\eta)^*=\overline{\sigma(\gamma,\eta)}\,b_\eta^* a_\gamma^*,
\qquad
b_\eta^*\star_\sigma a_\gamma^*=\sigma(-\eta,-\gamma)\,b_\eta^* a_\gamma^*,
\]
so $(a\star_\sigma b)^*=b^*\star_\sigma a^*$ for all homogeneous pairs if and only if \eqref{eq:cocycle-star-compat} holds,
and then for all finite sums by bilinearity.

(iii) The definition \eqref{eq:graded-product} uses only the homogeneous pieces and $\sigma$.
If two actions yield the same subspaces $A_\gamma$ (with the same labels), then the right-hand side of \eqref{eq:graded-product}
is identical for each homogeneous pair, hence the bilinear products coincide on the common graded core.
\end{proof}

\begin{proposition}[Recovering the action from the grading]\label{prop:B}
Let $K$ be compact abelian with $\Gamma=\widehat K$ and let $\alpha,\beta:K\to\mathrm{Aut}(A)$ be two actions such that
$A_{\mathrm{fin}}^\alpha$ and $A_{\mathrm{fin}}^\beta$ are defined and coincide as vector spaces.
Assume further that $A_{\mathrm{fin}}^\alpha$ is \emph{algebraically} the direct sum of its spectral subspaces and that
$\alpha$ and $\beta$ are determined by their values on $A_{\mathrm{fin}}^\alpha$ (e.g.\ $A_{\mathrm{fin}}^\alpha$ is dense in a chosen topology
and $\alpha,\beta$ are continuous).

\smallskip
\noindent (i) If $A_\gamma^\alpha=A_\gamma^\beta$ for all $\gamma\in\Gamma$ (with the same labeling), then $\alpha=\beta$.

\smallskip
\noindent (ii) Let $\psi\in\mathrm{Aut}(K)$ and write $\widehat\psi:\Gamma\to\Gamma$ for the induced automorphism
$\widehat\psi(\gamma):=\gamma\circ\psi^{-1}$.
Then $\beta=\alpha\circ\psi$ holds if and only if
\[
A_\gamma^\beta = A_{\widehat\psi(\gamma)}^\alpha\qquad(\gamma\in\Gamma),
\]
i.e.\ the gradings agree up to relabeling by $\widehat\psi$.
\end{proposition}

\begin{proof}
(i) Fix $t\in K$. For each $\gamma\in\Gamma$ and each $a\in A_\gamma^\alpha=A_\gamma^\beta$ we have
\[
\alpha_t(a)=\gamma(t)a=\beta_t(a).
\]
By linearity, $\alpha_t=\beta_t$ on the algebraic direct sum $A_{\mathrm{fin}}^\alpha$.
By the standing hypothesis that $\alpha$ and $\beta$ are determined by their values on $A_{\mathrm{fin}}^\alpha$
(e.g.\ by continuity and density), it follows that $\alpha_t=\beta_t$ on all of $A$ for every $t\in K$, hence $\alpha=\beta$.

(ii) Assume first that $\beta=\alpha\circ\psi$. Let $a\in A_\gamma^\beta$, so $\beta_t(a)=\gamma(t)a$ for all $t\in K$.
Then for all $t\in K$,
\[
\alpha_{\psi(t)}(a)=\beta_t(a)=\gamma(t)a.
\]
Set $s=\psi(t)$. Since $\psi$ is an automorphism, $t=\psi^{-1}(s)$, and the above becomes
\[
\alpha_s(a)=\gamma(\psi^{-1}(s))\,a=(\gamma\circ\psi^{-1})(s)\,a=\widehat\psi(\gamma)(s)\,a,
\]
so $a\in A_{\widehat\psi(\gamma)}^\alpha$. This shows $A_\gamma^\beta\subseteq A_{\widehat\psi(\gamma)}^\alpha$,
and the reverse inclusion follows by applying the same argument to $\psi^{-1}$.

Conversely, assume $A_\gamma^\beta=A_{\widehat\psi(\gamma)}^\alpha$ for all $\gamma$.
Take $a\in A_\gamma^\beta$. Then $a\in A_{\widehat\psi(\gamma)}^\alpha$, so for all $t\in K$,
\[
\alpha_{\psi(t)}(a)=\widehat\psi(\gamma)(\psi(t))\,a=\gamma(t)\,a=\beta_t(a).
\]
Thus $\beta_t=\alpha_{\psi(t)}$ on each spectral subspace and hence on $A_{\mathrm{fin}}^\beta$ by linearity.
By the determination hypothesis, $\beta=\alpha\circ\psi$ on $A$.
\end{proof}

\begin{corollary}[When two actions yield the same cocycle-twisted product]\label{cor:C}
Let $\sigma$ be a normalized unitary $2$-cocycle on $\Gamma=\widehat K$ and let $\alpha,\beta$ be as in Proposition~\ref{prop:B}.
\smallskip

\noindent (i) If $A_\gamma^\alpha=A_\gamma^\beta$ for all $\gamma$, then the cocycle twisted products on $A_{\mathrm{fin}}^\alpha=A_{\mathrm{fin}}^\beta$
defined by \eqref{eq:graded-product} coincide.

\smallskip
\noindent (ii) If $\beta=\alpha\circ\psi$ for $\psi\in\mathrm{Aut}(K)$, then the cocycle twisted product for $(\beta,\sigma)$ agrees with the cocycle twisted product for
\[
(\alpha,\sigma^\psi),\qquad
\sigma^\psi(\gamma,\eta):=\sigma\bigl(\widehat\psi^{-1}(\gamma),\widehat\psi^{-1}(\eta)\bigr).
\]
In particular, if $\sigma^\psi=\sigma$ (cocycle invariance under relabeling), then $(\alpha,\sigma)$ and $(\beta,\sigma)$ define the same product.
\end{corollary}

\begin{proof}
(i) If $A_\gamma^\alpha=A_\gamma^\beta$ for all $\gamma$, then $A_{\mathrm{fin}}^\alpha=A_{\mathrm{fin}}^\beta$ as graded vector spaces.
The cocycle twisted product is defined by the same homogeneous rule \eqref{eq:graded-product} on the same homogeneous pieces,
so the resulting bilinear products coincide.

(ii) Suppose $\beta=\alpha\circ\psi$. Let $a\in A_\gamma^\beta$ and $b\in A_\eta^\beta$.
By Proposition~\ref{prop:B}(ii), $a\in A_{\widehat\psi(\gamma)}^\alpha$ and $b\in A_{\widehat\psi(\eta)}^\alpha$.
Then the $(\alpha,\sigma^\psi)$-product gives
\[
a\star_{\sigma^\psi} b
=\sigma^\psi(\widehat\psi(\gamma),\widehat\psi(\eta))\,ab
=\sigma(\gamma,\eta)\,ab
=a\star_\sigma b
\]
(where the last product is the one defined using $\beta$ and $\sigma$).
Thus the products agree on homogeneous elements and hence on $A_{\mathrm{fin}}^\beta$ by bilinearity.
If $\sigma^\psi=\sigma$, this specializes to literal equality of the two cocycle-twisted products.
\end{proof}
\subsection{A $C^*$-level corollary: gradings are dual actions}\label{subsec:Cstar-grading}

\begin{theorem}[Topological gradings are implemented by the dual action]\label{thm:D}
Let $\Gamma$ be a discrete abelian group and let $A$ be a $C^*$-algebra equipped with closed subspaces $(A_\gamma)_{\gamma\in\Gamma}$ such that
\[
A_\gamma A_\eta\subseteq A_{\gamma+\eta},\qquad A_\gamma^*=A_{-\gamma},
\]
and the algebraic direct sum $A_{\mathrm{alg}}:=\bigoplus_{\gamma\in\Gamma}^{\mathrm{alg}}A_\gamma$ is dense in $A$.
Assume moreover that the canonical projection
\[
E_0:A_{\mathrm{alg}}\to A_0,\qquad E_0\Big(\sum_{\gamma}a_\gamma\Big):=a_0,
\]
extends by continuity to a faithful conditional expectation $E:A\to A_0$.

Then there exists a unique strongly continuous action $\widehat\alpha:\widehat\Gamma\to\mathrm{Aut}(A)$ such that
\[
A_\gamma=\{a\in A:\ \widehat\alpha_\chi(a)=\chi(\gamma)\,a\ \ \forall \chi\in\widehat\Gamma\}.
\]
Moreover, $E$ is given by Haar averaging:
\[
E(a)=\int_{\widehat\Gamma}\widehat\alpha_\chi(a)\,d\chi,\qquad a\in A.
\]
\end{theorem}
\begin{proof} This is standard in the theory of graded $C^*$-algebras and Fell bundles \cite{Exel}. Let $A_{\mathrm{alg}}:=\bigoplus_{\gamma\in\Gamma}^{\mathrm{alg}} A_\gamma$, which is a dense $*$-subalgebra of $A$ by hypothesis. For each character $\chi\in\widehat\Gamma$ define a linear map $\widehat\alpha_\chi:A_{\mathrm{alg}}\to A_{\mathrm{alg}}$ by
\begin{equation}\label{eq:gauge-on-alg}
\widehat\alpha_\chi\!\left(\sum_{\gamma} a_\gamma\right):=\sum_{\gamma} \chi(\gamma)\,a_\gamma,
\qquad a_\gamma\in A_\gamma,
\end{equation}
where the sum is finite. This is well-defined because elements of $A_{\mathrm{alg}}$ have unique homogeneous expansions.
A direct computation on homogeneous elements shows that $\widehat\alpha_\chi$ is a $*$-homomorphism:
if $a_\gamma\in A_\gamma$ and $b_\eta\in A_\eta$, then
\[
\widehat\alpha_\chi(a_\gamma b_\eta)=\chi(\gamma+\eta)a_\gamma b_\eta
=\chi(\gamma)\chi(\eta)a_\gamma b_\eta=\widehat\alpha_\chi(a_\gamma)\widehat\alpha_\chi(b_\eta),
\]
and similarly $\widehat\alpha_\chi(a_\gamma^*)=\overline{\chi(\gamma)}\,a_\gamma^*=\chi(-\gamma)\,a_\gamma^*$.
Moreover $\widehat\alpha_{\chi^{-1}}$ is the inverse of $\widehat\alpha_\chi$ on $A_{\mathrm{alg}}$.

Consider the right Hilbert $A_0$-module $X_0:=A_{\mathrm{alg}}$ with $A_0$-valued inner product
\[
\langle x,y\rangle_{A_0}:=E(x^*y)\in A_0,
\]
and let $X$ be its completion. Faithfulness of $E$ implies $\langle x,x\rangle_{A_0}=0$ only when $x=0$ in $A$,
so this is a genuine inner product.
Left multiplication defines a $*$-homomorphism $\lambda:A\to\mathcal{L}(X)$ by $\lambda(a)(x)=ax$ on $X_0$.
We claim $\lambda$ is faithful (hence isometric).
Indeed, suppose $\lambda(a)=0$. Then $ax=0$ in $X$ for all $x\in X_0$, so
\[
0=\langle ax,ax\rangle_{A_0}=E\big((ax)^*(ax)\big)=E(x^*a^*a x)\qquad(x\in X_0).
\]
By faithfulness of $E$, $E(b^*b)=0\Rightarrow b=0$, hence $ax=0$ in $A$ for all $x\in A_{\mathrm{alg}}$.
Since $A_{\mathrm{alg}}$ is dense in $A$ and the map $x\mapsto ax$ is norm-continuous, we get $ax=0$ for all $x\in A$.
Let $(e_i)$ be a contractive approximate identity of $A$. Then $a=\lim_i a e_i=0$.
Thus $\ker\lambda=\{0\}$, as claimed.
Now define $U_\chi:X_0\to X_0$ by the same formula as \eqref{eq:gauge-on-alg}.
Since $E$ extends $E_0$, we have $E(A_\delta)=0$ for $\delta\neq 0$; using also that $A_\gamma^*A_\eta\subseteq A_{\eta-\gamma}$, one checks
\[
\langle U_\chi x, U_\chi y\rangle_{A_0}
=E\!\left((U_\chi x)^*(U_\chi y)\right)
=E(x^*y)
=\langle x,y\rangle_{A_0}.
\]
Thus $U_\chi$ extends to a unitary on $X$.
Moreover, for $a\in A_{\mathrm{alg}}$ one has
\[
\lambda(\widehat\alpha_\chi(a))=U_\chi\,\lambda(a)\,U_\chi^*.
\]
Hence, for $a\in A_{\mathrm{alg}}$,
\[
\|\widehat\alpha_\chi(a)\|
=\|\lambda(\widehat\alpha_\chi(a))\|
=\|U_\chi\lambda(a)U_\chi^*\|
=\|\lambda(a)\|
=\|a\|.
\]
So each $\widehat\alpha_\chi$ is isometric on $A_{\mathrm{alg}}$ and therefore extends uniquely by continuity to a
$*$-automorphism of $A$, still denoted $\widehat\alpha_\chi$.

For $a\in A_{\mathrm{alg}}$, the map $\chi\mapsto \widehat\alpha_\chi(a)$ is continuous because it is a finite sum of continuous characters.
For general $a\in A$, approximate by $a^{(n)}\in A_{\mathrm{alg}}$ and use isometry:
\[
\|\widehat\alpha_\chi(a)-\widehat\alpha_{\chi_0}(a)\|
\le
\|a-a^{(n)}\|+\|\widehat\alpha_\chi(a^{(n)})-\widehat\alpha_{\chi_0}(a^{(n)})\|+\|a^{(n)}-a\|,
\]
so strong continuity follows.
The group law $\widehat\alpha_{\chi\chi'}=\widehat\alpha_\chi\circ\widehat\alpha_{\chi'}$ holds on $A_{\mathrm{alg}}$ by inspection
and hence on $A$ by continuity.

If $a\in A_\gamma$, then \eqref{eq:gauge-on-alg} gives $\widehat\alpha_\chi(a)=\chi(\gamma)a$ for all $\chi$, hence
$A_\gamma\subseteq \{a\in A:\widehat\alpha_\chi(a)=\chi(\gamma)a\}$.
Conversely, define the Fourier coefficient map
\[
P_\gamma(a):=\int_{\widehat\Gamma}\overline{\chi(\gamma)}\,\widehat\alpha_\chi(a)\,d\chi,
\]
where $d\chi$ is normalized Haar measure and the integral is a Bochner integral in $A$.
Since $\widehat\alpha_\chi$ is isometric and strongly continuous, $P_\gamma$ is a bounded linear map with $\|P_\gamma\|\le 1$.
On $A_{\mathrm{alg}}$ one checks directly that $P_\gamma(\sum_\eta a_\eta)=a_\gamma$, hence $P_\gamma(A)\subseteq A_\gamma$ by density
and closedness of $A_\gamma$. If $\widehat\alpha_\chi(a)=\chi(\gamma)a$ for all $\chi$, then $P_\gamma(a)=a$, hence $a\in A_\gamma$.
This proves the spectral subspace characterization.

Finally, for $a\in A_{\mathrm{alg}}$,
\[
\int_{\widehat\Gamma}\widehat\alpha_\chi(a)\,d\chi = a_0,
\]
because $\int_{\widehat\Gamma}\chi(\gamma)\,d\chi=0$ for $\gamma\neq 0$ and equals $1$ for $\gamma=0$.
Since $E$ extends $E_0$, we have $E(\sum_\gamma a_\gamma)=a_0$ for $\sum_\gamma a_\gamma\in A_{\mathrm{alg}}$, and hence
$E(a)=\int_{\widehat\Gamma}\widehat\alpha_\chi(a)\,d\chi$ for all $a\in A$ by continuity.

Uniqueness of $\widehat\alpha$ follows because it is determined on the dense graded subalgebra by the eigenvalue relations
$\widehat\alpha_\chi|_{A_\gamma}=\chi(\gamma)\mathrm{id}$.
\end{proof}

\begin{corollary}[No genuinely action-free graded cocycle twists at the $C^*$ level]\label{cor:E}
In the setting of Theorem~\ref{thm:D}, any graded cocycle twist determined by a $2$-cocycle $\sigma$ on $\Gamma$
is (on the dense graded $*$-subalgebra $\bigoplus_\gamma^{\mathrm{alg}}A_\gamma$) given by the homogeneous rule
$a_\gamma\star_\sigma b_\eta=\sigma(\gamma,\eta)a_\gamma b_\eta$ and is precisely the cocycle deformation associated to the dual action
$\widehat\alpha$.

In particular, if $A=C(X)$ is commutative and admits a nontrivial topological $\Gamma$-grading with a noncommutative cocycle twist,
then $X$ admits a nontrivial continuous action of the compact abelian group $\widehat\Gamma$ (via Gelfand duality),
so noncommutativity in the graded cocycle class forces the presence of compact abelian symmetry at the $C^*$-level.
\end{corollary}

\begin{proof}
Let $A_{\mathrm{alg}}=\bigoplus_{\gamma}^{\mathrm{alg}}A_\gamma$ be the dense graded $*$-subalgebra.
For a normalized $2$-cocycle $\sigma$ on $\Gamma$, the graded cocycle twisted product on $A_{\mathrm{alg}}$ is defined by the homogeneous rule
$a_\gamma\star_\sigma b_\eta=\sigma(\gamma,\eta)a_\gamma b_\eta$.

By Theorem~\ref{thm:D}, the topological $\Gamma$-grading is implemented by a strongly continuous action
$\widehat\alpha:\widehat\Gamma\to\mathrm{Aut}(A)$ whose spectral subspaces are precisely the $A_\gamma$.
In the standard cocycle deformation of a $C^*$-dynamical system (in particular, in the cocycle/twist formalism
appearing in Kasprzak-type deformations), the deformed multiplication on spectral subspaces is given by exactly the same cocycle factor.
Thus, on the canonical dense graded core $A_{\mathrm{alg}}$, the cocycle deformation associated to $\widehat\alpha$
coincides with the graded cocycle twisted product above.

For the commutative case $A=C(X)$, any $C^*$-action $\widehat\alpha$ of the compact abelian group $\widehat\Gamma$
corresponds by Gelfand duality to a continuous action of $\widehat\Gamma$ on $X$.
If the $\Gamma$-grading is nontrivial and the cocycle twist is noncommutative, then there exist nonzero homogeneous components
in at least two degrees and the implementing action $\widehat\alpha$ is necessarily nontrivial; hence $X$ admits a nontrivial
continuous $\widehat\Gamma$-action.
\end{proof}

\section{Outlook}\label{sec:outlook}

This paper develops a symmetry-free spectral channel-twist framework, proves unconditional well-posedness on the finite spectral core,
and cleanly separates analytic closure from algebraic constraints via an explicit Sobolev hypothesis.  It also clarifies the structural role
of grading in classical strict deformation pictures: in situations where the relevant abelian action has a discrete spectral decomposition (so that the algebra admits a character-indexed homogeneous decomposition), cocycle twisting is governed by
the induced grading and can be realized as a refined instance of our channel twist. The spectral-channel deformation is always defined on the finite spectral core and in $L^2(M)$, independently of any homogeneous decomposition. However, the graded cocycle classification and the uniform recovery of classical strict deformation products (Rieffel, Connes-Landi and Kasprzak) rely on the existence of a discrete spectral decomposition of a compatible abelian group action.

The main directions suggested by the present results are the following.

\medskip
\noindent\textbf{(1) Analytic criteria for Sobolev closure.}
Assumption~\ref{ass:A} isolates the analytic input needed to iterate $\star_\omega$ on a Sobolev algebra and to interpret associativity
as an identity in a closed function space.  Developing geometric/spectral conditions that imply such boundedness (for instance, via bilinear
spectral multiplier estimates adapted to the channel decomposition) would turn the current associativity criterion into a practical existence tool.

\medskip
\noindent\textbf{(2) Beyond scalar phases: multiplicity and tensor-category deformations.}
The gauge and obstruction results indicate strong rigidity in the scalar regime, and the graded cocycle corner largely overlaps with
classical cocycle twisting mechanisms in discrete spectral decomposition settings.  To produce genuinely new intrinsic noncommutative strictly associative examples,
one is naturally led to non-scalar twisting on higher-dimensional channel/multiplicity spaces, where associativity is governed by nontrivial
coherence constraints.  A forthcoming paper develops such a deformation framework organized by rigid $C^*$-tensor categories and
multiplicity/intertwiner data, designed to subsume the graded cocycle corner as a special case and to enable new example construction.

\bibliographystyle{unsrt}
\bibliography{refs}

\end{document}